\documentclass[11pt,a4paper,oneside,draft]{article}

\usepackage{a4,amsthm,amsfonts}
\usepackage[intlimits]{amsmath}

\author{Stefan Boller\\
\\
Graduiertenkolleg Quantenfeldtheorie\\
Universit\"at Leipzig, Augustusplatz 10/11\\
04109 Leipzig, Germany\\
e-mail: Stefan.Boller@itp.uni-leipzig.de}

\title{Characterization of Cyclic and Separating Vectors and Application to an
  Inverse Problem in Modular Theory\\
  I. Finite Factors}

\newtheorem{thm}{Theorem}[section]
\newtheorem{cor}[thm]{Corollary}
\newtheorem{lem}[thm]{Lemma}
\newtheorem{prop}[thm]{Proposition}

\theoremstyle{definition}
\newtheorem{defn}{Definition}[section]

\theoremstyle{remark}
\newtheorem{rem}{Remark}[section]

\newtheorem{bsp}{Example}[section]

\numberwithin{equation}{section}

\newcommand{\thmref}[1]{Theorem~\ref{#1}}
\newcommand{\secref}[1]{\S\ref{#1}}
\newcommand{\lemref}[1]{Lemma~\ref{#1}}
\newcommand{\corref}[1]{Corollary~\ref{#1}}
\newcommand{\propref}[1]{Proposition~\ref{#1}}
\newcommand{\defnref}[1]{Definition~\ref{#1}}
\newcommand{\remref}[1]{Remark~\ref{#1}}

\newcommand{\bspref}[1]{Example~\ref{#1}}

\DeclareMathOperator{\ad}{ad}
\DeclareMathOperator{\aut}{aut}

\DeclareMathOperator{\tr}{tr}

\newcommand{\interval}[1]{\mathinner{#1}}

\newcommand{\abge}[2]{\interval{\lbrack#1,#2\rbrack}}
\newcommand{\offenl}[2]{\interval{(#1,#2\rbrack}}
\newcommand{\offenr}[2]{\interval{\lbrack#1,#2)}}

\newcommand{\eval}[2][\right]{\relax
  \ifx#1\right\relax \left.\fi#2#1\rvert}

\let\abs=\envert

\let\norm=\enVert

\newcommand{\spitz}[1]{\bigl\langle#1\bigr\rangle}
\newcommand{\SProd}[2]{\spitz{#1\rvert#2}}

\newcommand{\field}[1]{\mathbb{#1}}
\newcommand{\R}{\field{R}}
\newcommand{\C}{\field{C}}

\newcommand{\N}{\field{N}}

\newcommand{\kon}[1]{\overline{#1}}

\newcommand{\clos}[1]{\overline{#1}}

\newcommand{\raum}[1]{\mathcal{#1}}
\renewcommand{\H}{\raum{H}}

\newcommand{\alg}[1]{\mathcal{#1}}
\newcommand{\M}{\alg{M}}

\newcommand{\Op}[1]{\mathrm{#1}}

\newcommand{\Dom}[1]{\mathcal{D}(\Op{#1})}

\newcommand{\Ran}[1]{\mathcal{W}(\Op{#1})}

\newcommand{\Lim}[1]{\lim_{#1\rightarrow\infty}}

\begin{document}

  \maketitle

  \begin{abstract}  
    In this paper we examine an inverse problem in the modular theory of von
    Neumann algebras in the case
    of finite factors. First we give a characterization of cyclic and
    separating vectors for finite factors in terms of operators
    associated with this vector and being
    affiliated with the factor. Further we show how this operator generates
    the 
    modular objects of the given cyclic and separating vector generalizing an
    idea of Kadison and Ringrose. With the help
    of these rather technical results we show under an appropriate condition,
    which is always fulfilled for finite type I factors, that 
    there exists another simple class of solutions of the inverse problem
    beside a trivial one which always exists. Finally we give a complete
    classification of the solutions of the inverse problem in the case of
    modular operators having pure point
    spectrum which is no restriction in the type I case. In a subsequent paper
    these results will be generalized to all semifinite factors.
  \end{abstract}

  \setcounter{page}{1}

\section{The Inverse Problem in Modular Theory}
\label{sec:1}

Let $\M_{0}$ be a von Neumann algebra on a separable Hilbert space $\H_{0}$
with a cyclic and separating vector $u_{0}$. Then modular theory shows the
existence of a modular operator $\Delta_{0}$ and a modular conjugation
$\Op{J}_{0}$ (the modular objects $(\Delta_{0},\Op{J}_{0})$) belonging to the
vector $u_{0}$. In this paper we examine the inverse problem of constructing
algebras $\M$ having the same cyclic and separating vector and modular objects
as $\M_{0}$:

\textbf{The Inverse Problem}

Let $(\Delta_{0},\Op{J}_{0})$ be the modular objects for the von Neumann
algebra $\M_{0}$ with cyclic and separating vector $u_{0}$. Characterize all
von Neumann algebras $\M$ isomorphic to $\M_{0}$ with the following
properties:
\begin{enumerate}
  \item $u_{0}$ is also cyclic and separating for $\M$,
  \item $(\Delta_{0},\Op{J}_{0})$ are the modular objects for $(\M,u_{0})$.
\end{enumerate}
Let $NF_{\M_{0}}(\Delta_{0},\Op{J}_{0},u_{0})$ denote all solutions $\M$ of
the inverse problem.

We consider here only von Neumann factors instead of general von Neumann
algebras. 

General facts of this problem were considered by Wollenberg (s. \cite{WolI},
\cite{WolII}) in the frame of algebraic quantum field theory (cf. also
\cite{WolIII}, \cite{BOR}). In
\cite{WolI} the following useful proposition was proven:
\begin{prop}\label{prop1:1}
  A von Neumann factor $\M$ belongs to
  $NF_{\M_{0}}(\Delta_{0},\Op{J}_{0},u_{0})$ iff there is a unitary operator
  $\Op{U}$ s.t.
  \begin{enumerate}
    \item $\Op{M}=\Op{U}\M_{0}\Op{U}^{*}$,
    \item $u:=\Op{U}^{*}u_{0}$ is a cyclic and separating vector for $\M_{0}$,
    \item $(\Op{U}^{*}\Delta_{0}\Op{U},\Op{J}_{0})$ are the modular objects
      for $(\M_{0},u)$,
    \item $\Op{U}$ commutes with $\Op{J}_{0}$.
  \end{enumerate}
\end{prop}
Further Wollenberg showed in \cite{WolII} that for finite type $I$ factors there
exists always two simple classes of solutions of the inverse problems. Further
he showed that in this case the solution of the inverse problem can be reduced to an
algebraic problem. 

In this paper we consider finite factors and hence generalize some of those
results for type $II_{1}$ factors and slightly modify the representation of
the results. In particular we establish a correspondence between cyclic and
separating vectors for a finite von Neumann factor and invertible operators
affiliated with this factor. Then we show how the modular operator can be
constructed with the help of this operator generalizing an idea of Kadison and
Ringrose (\cite[9.6.11]{KRII}). In the forth section this result will be used
to show that the second of Wollenberg's simple classes exists in some cases
also for type $II_{1}$ factors, but not always. In the fifth section we give a
complete classification of the solutions of the inverse problem when the
generating operator has pure point spectrum.

Notice that in this paper all Hilbert spaces are separable, i.e. the von
Neumann algebras are countably decomposable.

\section{Characterization of Vectors by Affiliated Operators}
\label{sec:2}

Throughout this section $(\M_{0},\H_{0})$ is a finite von Neumann factor
with a cyclic, separating vector $u_{0}\in\H_{0}$ and (unique) tracial state
$\tr$ continued to all the positive operators affiliated with $\M_{0}$. In this
section we establish a connection between operators affiliated
with a factor and vectors of the underlying Hilbert space. These operators are
characterized by the finite trace of the square of their modulus. Further a
vector is cyclic iff the corresponding operator is injective, and it is
separating iff the corresponding operator has dense range. Last we will show
that every operator with the finite trace condition gives rise to a vector
with the corresponding properties.

Since $\M_{0}$ has a cyclic and separating vector, the tracial
state $\tr$ is a vector state generated by a cyclic (and hence
separating) trace vector $u_{\tr}$.  Using this vector the trace can be
continued to all the positive closed operators $\Op{A}$ affiliated  with
$\M_{0}$ by 
\begin{equation}\label{eq2:1:1}
  \tr(\Op{A}):=\int\lambda d\norm{\Op{E}_{\lambda}u_{\tr}}^{2},
\end{equation}
where $\Op{E}_{\lambda}$ is the spectral measure of $\Op{A}$.

It is known [cf. \cite[8.7.60(vi)]{KRII} in the
case of finite algebras that we can reach every vector from a cyclic vector by
an operator affiliated with $\M_{0}$. In particular we can give the following
definition:
\begin{defn}\label{defn2:1:1}
  For every vector $u\in\H_{0}$ we denote with $\Op{T}_{u}$ an operator
  affiliated with $\M_{0}$ s.t. $u_{\tr}\in\Dom{T_{u}}$ and
  $\Op{T}_{u}u_{\tr}=u$.
\end{defn}

Later we will show that the operator defined in \defnref{defn2:1:1} is unique.
The next proposition shows some useful properties of the domain of operators
affiliated with a von Neumann algebra: 
\begin{prop}\label{prop2:1:2}
  Let $\M_{0}$ be a von Neumann algebra and $\Op{T}\eta\M_{0}$.
  \begin{enumerate}
    \item Let $u_{\tr}\in\H_0$ be a trace vector in $\H_{0}$ and
      $u_{\tr}\in\Dom{T}$, then
      \begin{enumerate}
      \item 
        $u_{\tr}\in\Dom{T^{*}}$ and
        $u_{\tr}\in\Dom{(\Op{T}^{*}\Op{T})^{1/2}}$,
      \item
        $\M_{0}u_{\tr}\subset\Dom{T}$, $\M_{0}u_{\tr}\subset\Dom{T^{*}}$, and
        $\M_{0}u_{\tr}\subset\Dom{(\Op{T}^{*}\Op{T})^{1/2}}$.
      \item 
        If $u_{\tr}$ is also cyclic then $\M_{0}u_{\tr}$ is a core for
        $\Op{T}$, $\Op{T}^{*}$, and $(\Op{T}^{*}\Op{T})^{1/2}$, and the
        assertions of b) and c) hold also for $\M_{0}^{'}$ instead of $\M_{0}$.
      \end{enumerate}
    \item Let $u\in\Dom{T}\subset\H_{0}$. Then 
      \begin{enumerate}
      \item
        $\M_{0}^{'}u\subset\Dom{T}$ and
        $\M_{0}^{'}u\subset\Dom{(\Op{T}^{*}\Op{T})^{1/2}}$, and
      \item       
        $[\M_{0}^{'}\Op{T}u]=[\Op{T}\M_{0}^{'}u]$,
         and
        $[\M_{0}^{'}(\Op{T}^{*}\Op{T})^{1/2}u]=
        [(\Op{T}^{*}\Op{T})^{1/2}\M_{0}^{'}u]$.  
    \end{enumerate}
  \end{enumerate}
\end{prop}

\begin{proof}
  \begin{enumerate}
 \item
    Let $\Op{T}=\Op{VH}$ the polar decomposition of $\Op{T}$ and
    $\Op{E}_{\lambda}\in\M_{0}$ the spectral resolution of $\Op{H}$. Then
    $\Dom{T}=\Dom{H}=\Dom{(\Op{T}^{*}\Op{T})^{1/2}}$ and
    $u_{\tr}\in\Dom{H}$, i.e.
    \begin{equation*}
      \int\lambda^{2}d\norm{\Op{E}_{\lambda}u_{\tr}}^{2}<\infty.
    \end{equation*}
    Set now $u:=\Op{U}u_{\tr}$ with a unitary $\Op{U}\in\M_{0}$. Since
    $u_{\tr}$ is a trace vector we have
    \begin{equation*}
      \begin{split}
        \int\lambda^{2}d\norm{\Op{E}_{\lambda}u}^{2}
        &=\int\lambda^{2}d\norm{\Op{E}_{\lambda}\Op{U}u_{\tr}}^{2}\\
        &=\int\lambda^{2}
        d\norm{\Op{U}^{*}\Op{E}_{\lambda}u_{\tr}}\\
        &=\int\lambda^{2}d\norm{\Op{E}_{\lambda}u_{\tr}}^{2}<\infty,
      \end{split}
    \end{equation*}
    i.e. $u\in\Dom{H}$. Since every element of $\M_{0}$ is the linear
    combination of at
    most 4 unitaries, it follows, that
    $\M_{0}u_{\tr}\subset\Dom{H}=\Dom{(\Op{T}^{*}\Op{T})^{1/2}}=\Dom{T}$.

    Since $\Op{T}^{*}=\Op{H}\Op{V}^{*}$, $u_{\tr}\in\Dom{H}$, and
    $\Op{V}^{*}\in\M_{0}$, the assertions follow also for $\Op{T}^{*}$
    applying the results just proven.  

    Let now $u_{\tr}$ be cyclic for $\M_{0}$. Then for every
    $\Op{M}^{'}\in\M_{0}^{'}$ there exists exactly one $\Op{M}\in\M_{0}$
    s.t. $\Op{M}^{'}u_{\tr}=\Op{M}u_{\tr}$, hence
    $\M_{0}^{'}u_{\tr}=\M_{0}u_{\tr}$. Further $\M_{0}u_{\tr}$ is dense in
    $\H_{0}$ and $\M_{0}u_{\tr}\subset\Dom{T}$, hence
    $\clos{\Op{T}_{\M_{0}u_{\tr}}}\subset\Op{T}$. Since $\M_{0}$ is finite (it
    possesses a cyclic trace vector) it follows from this
    \begin{equation*}
      \clos{\Op{T}_{\M_{0}u_{\tr}}}=\Op{T},
    \end{equation*}
    i.e. $\M_{0}u_{\tr}$ (and $\M_{0}^{'}u_{\tr}$) is a core for $\Op{T}$ and,
    similarly, also for $(\Op{T}^{*}\Op{T})^{1/2}$ and $\Op{T}^{*}$. 
  \item
    This first part follows in the same way as the corresponding assertions in
    1.
 
    Now $\Op{M}^{'}u\in\Dom{T}$ for every
    $\Op{M}^{'}\in\M_{0}^{'}$ and, since $\Op{T}\eta\M_{0}$,
    $\Op{T}\Op{M}^{'}u=\Op{M}^{'}\Op{T}u$
    s.t. $[\M_{0}^{'}\Op{T}u]=[\Op{T}\M_{0}^{'}u]$.
  \end{enumerate}
\end{proof}

The next proposition shows that separating vectors are also separating for
operators only affiliated with a von Neumann algebra:
\begin{prop}\label{prop2:1:1}
  Let $u\in\H_0$ be a separating vector for a von Neumann algebra $\M_0$. Let
  further $\Op{A}\eta\M_0$, $u_{0}\in\Dom{A}$ and $\Op{A}u_0=0$. Then
  $\Op{A}=0$.
\end{prop}

\begin{proof}
  Let $\Op{A}=\Op{VH}$ with $\Op{V}\in\M_0$ be a partial isometry and
  $0\leq\Op{H}\eta\M_0$, the polar decomposition of $\Op{A}$. Then we have
  \begin{equation*}
    0=\Op{A}u_0=\Op{VH}u_0\quad\Rightarrow\quad\Op{H}u_0=0,
  \end{equation*}
  since $\Op{V}$ is a partial isometry from $\clos{\Ran{H}}$ to
  $\clos{\Ran{A}}$. Let now $\Op{E}_\lambda$ be the spectral measure for
  $\Op{H}$. Then
  \begin{equation*}
    0=\Op{H}u_0=\int_{\R_\ge}\lambda d\Op{E}_\lambda u_0.
  \end{equation*}
  From this we see, that $\Op{E}_{\R_>}u_0=0$ and, since $u_0$ is separating
  for $\M_0$ and $\Op{E}_{\R_>}\in\M_0$, $\Op{E}_{\R_>}=0$. Since $\Op{H}$ is
  positive $\Op{E}_{\R_<}=0$,too, and therefore $\Op{E}_{\{0\}}=\Op{I}$ and
  $\Op{H}=0$.
\end{proof}

Now we can show the uniqueness of $\Op{T}_{u}$:
\begin{cor}\label{cor2:1:1}
  The operator $\Op{T}_{u}$ defined in \defnref{defn2:1:1} is unique.
\end{cor}
\begin{proof}
  Let $\Op{T}$ and $\Op{S}$ be two operators affiliated with $\M_{0}$, s.t.
  \begin{equation*}
    \Op{T}u_{\tr}=\Op{S}u_{\tr}=u.
  \end{equation*}
  Then $\Op{T}-\Op{S}$ is closable, since $\M_{0}$ is finite (cf.
  \cite[8.7.60]{KRII}), and its closure $\clos{\Op{T}-\Op{S}}$ is affiliated
  with $\M_{0}$. Then it follows from 
  \propref{prop2:1:1} that $\clos{\Op{T}-\Op{S}}=0$, i.e. $\Op{T}$ and
  $\Op{S}$ agree on the intersection of their domains. Since $\M_{0}u_{\tr}$
  is a core for both $\Op{S}$ and $\Op{T}$ (cf. \propref{prop2:1:2}), they are
  equal.
\end{proof}

For the proof of the main lemma in this subsection the next proposition is
also necessary:
\begin{prop}\label{prop2:1:3}
Let $\M_0$ be a von Neumann algebra on the Hilbert space $\H_0$ and
$\Op{T}\eta\M_0$. Then
\begin{equation}\label{eq2:1:6}
  [\M_0\Op{T}u]\subset[\M_0u]\quad\forall u\in\Dom{T}\subset\H_0
\end{equation}
\end{prop}

\begin{proof}
  Let $\Op{T}=\Op{VH}$ be the polar decomposition of $\Op{T}$ and 
  $\Op{E}_n\in\M_0$ the spectral projection of $\Op{H}$ for the interval
  $\abge{-n}{n}$. Then $\Op{TE}_n=\Op{VHE}_n\in\M_0$, s.t.  
  \begin{equation*}
    \M_0\Op{TE}_n\subset\M_0\text{ and }
    \M_0\Op{TE}_n u\subset\M_0 u\quad\forall u\in\H_0.
  \end{equation*}
  Now 
  \begin{equation*}
    \Op{MT}u=\Lim{n}\Op{MTE}_n u\in[\M_0 u]\quad\forall u\in\Dom{T},
  \end{equation*}
  since $\bigcup_{n=1}^{\infty}\Op{E}_{n}(\H)$ is a core for $\Op{H}$ and
  $\Op{MTE}_n u\in\M_0\Op{TE}_n u\subset\M_0 u$. 
\end{proof}

Now we can formulate and proof the main statement of this subsection:
\begin{lem}\label{lem2:1:1}
  Let $\Op{T}_{u}$ be the operator defined in \defnref{defn2:1:1}. Then
  \begin{enumerate}
    \item $\tr(\Op{T}_{u}\Op{T}_{u}^{*})=\tr(\Op{T}_{u}^{*}\Op{T}_{u})<\infty$.
    \item $u$ is cyclic iff $\Op{T}_{u}$ is injective.
    \item $u$ is separating iff $\Op{T}_{u}$ has dense range.
    \item $u$ is cyclic and separating iff $\Op{T}_{u}$ is injective and has
  dense range, i.e. iff $\Op{T}_{u}$ is invertible. 
  \end{enumerate}
\end{lem}

\begin{proof}
  \begin{enumerate}
    \item
      Since $u_{\tr}\in\Dom{T_{u}}$ also $u_{\tr}\in\Dom{H}$,
      where $\Op{T}_{u}=\Op{V}\Op{H}$ is the polar decomposition of
      $\Op{T}_{u}$. This is equivalent with
      \begin{equation*}
        \begin{split}
          \infty
          &>\int\lambda^{2} d\norm{\Op{E}_{\lambda}u_{\tr}}^{2}
          =\tr(\Op{H}^{2})=\tr(\Op{T}_{u}^{*}\Op{T}_{u})\\
          &=\int\lambda^{2} d\norm{\Op{E}_{\lambda}\Op{V}^{*}u_{\tr}}^{2}=
          \tr(\Op{V}\Op{H}\Op{H}\Op{V}^{*})=\tr(\Op{T}_{u}\Op{T}_{u}^{*}),
        \end{split}
      \end{equation*}
      where $\Op{E}_{\lambda}$ is the spectral measure of $\Op{H}$ and we
      have used the trace property of $u_{\tr}$ and that
      $\Op{E}_{\lambda}u_{\tr}$ is in the
      final space of $\Op{V}$.
     \item If $u$ is cyclic there is an operator $\Op{S}\eta\M_{0}$,
      s.t. $u\in\Dom{S}$ and 
      \begin{equation*}
        \Op{S}u=u_{\tr}.
      \end{equation*}
      Now we have the following equality
      \begin{equation}\label{eq2:1:3}
        \Op{S}\cdot\Op{T}_{u}u_{\tr}=\Op{S}u=u_{\tr}\quad\Leftrightarrow\quad
        (\Op{S}\cdot\Op{T}_{u}-\Op{I})u_{\tr}=0
      \end{equation}
      Since $\M_{0}$ is finite $\Op{S}\cdot\Op{T}_{u}$
      is densely defined and closable
      and the closure $\Op{S}\Op{T}_{u}\eta\M_0$
      (cf.\cite[8.7.60(iii)]{KRII}). Furthermore also
      $\Op{S}\Op{T}_{u}-\Op{I}$ is closable and its closure $\Op{A}$ is
      affiliated with $\M_{0}$. From this, \eqref{eq2:1:3} and
      \propref{prop2:1:1} it follows $\Op{A}=0$ since $u_{\tr}$ is separating,
      and therefore $\Op{T}_{u}$ is injective. 

      Suppose now $\Op{T}_{u}$ is injective, i.e. $\Op{H}$ and $\Op{H}^{-1}$
      are both affiliated with $\M_{0}$, where $\Op{T}_{u}=\Op{V}\Op{H}$ is
      the polar decomposition of $\Op{T}_{u}$. Now [cf. \propref{prop2:1:3}] 
      \begin{equation*}
            [\M_0\Op{H}^{-1}x]\subset[\M_0x]\quad\forall x\in\Dom{H^{-1}}.
       \end{equation*}
       So we have the following chain
       \begin{equation*}
         \H_{0}=[\M_0u_{\tr}]=[\M_0\Op{H}^{-1}\Op{H}u_{\tr}]
         \subset[\M_0\Op{H}u_{\tr}]=[\M_0\Op{V}\Op{H}u_{\tr}]\subset\H_{0},
       \end{equation*}
       since $u_{\tr}$ is cyclic for $\M_0$ and $\Op{V}$ is a partial isometry
       with $\Op{H}u_{\tr}$ in its initial space.
       This means that $u=\Op{T}_{u}u_{\tr}$ is cyclic for $\M_0$.
   \item If $u$ is separating for $\M_{0}$ it is cyclic for the commutant
     $\M_{0}^{'}$, i.e.
     \begin{equation*}
       \H_{0}=[\M_{0}^{'}u]=[\M_{0}^{'}\Op{T}_{u}u_{\tr}]
       =[\Op{T}_{u}\M_{0}^{'}u_{\tr}].
     \end{equation*}
     The last equality follows from \propref{prop2:1:2}. This shows that
     $\Op{T}_{u}$ has dense range. 

     Suppose now that $\Op{T}_{u}$ has dense range.
       Let $\Op{A}\in\M_0$ and $\Op{A}u=0$. This means
       $\Op{A}(\Op{T}_{u}u_{\tr})=0$. Since $\Op{T}_{u}\eta\M_0$ and $\M_0$ is
       finite we know that $\Op{AT}_{u}\eta\M_0$. And using
       \propref{prop2:1:1} we derive $\Op{A}\cdot\Op{T}_{u}=0_{\Dom{T_{u}}}$
       and, since $\Op{T}_{u}$ has dense range and $\Op{A}$ is bounded,
       $\Op{A}=0$, s.t. $u=\Op{T}_{u}u_{\tr}$ is separating.
   \item This follows from 1. and 2.
 \end{enumerate}
\end{proof}

\begin{rem}\label{rem2:1:1}
  \begin{enumerate}
  \item The finite trace condition of \lemref{lem2:1:1} is not only necessary
    for an operator being the operator corresponding to a vector in the sense
    of \defnref{defn2:1:1} but also sufficient, how the following short
    calculation shows: Let $\Op{T}\eta\M_{0}$ and
    $\tr(\Op{T}^{*}\Op{T})<\infty$, i.e.
    \begin{equation*}
      \infty>\tr(\Op{T}^{*}\Op{T})=\tr(\Op{H}^{2})=
      \int\lambda^{2}d\norm{\Op{E}_{\lambda}u_{\tr}}^{2}, 
    \end{equation*}
    where $\Op{T}=\Op{V}\Op{H}$ is the polar decomposition of $\Op{T}$ and
    $\Op{E}_{\lambda}$ the spectral resolution of $\Op{H}$. This shows, that
    $u_{\tr}\in\Dom{H}=\Dom{T}$. Now $u:=\Op{T}u_{\tr}$ and \corref{cor2:1:1}
    shows that $\Op{T}$ is the
    unique operator associated with $u$ in the sense of \defnref{defn2:1:1}.
  \item If $L_{2}(\M_{0},\tr)$ is the Hilbert space obtained from $\M_{0}$ by
    completion w.r.t. the trace norm, then the finite trace condition of
    \lemref{lem2:1:1} is equivalent to $\Op{T}_{u}\in L_{2}(\M_{0},\tr)$,
    i.e. $\Op{T}_{u}$ is a quadratic integrable operator affiliated with
    $\M_{0}$. 
  \end{enumerate}
\end{rem}

The results of this section can be subsumed in the next two theorems:
\begin{thm}\label{thm2:1}
  Let ($\M_{0},\H_{0}$) be a finite von Neumann factor. Let further
  $u\in\H_{0}$. Then there is
  exactly one operator $\Op{T}_{u}\eta\M_{0}$ associated with the vector $u$
  in the sense of \defnref{defn2:1:1}, having the
  following properties: 
  \begin{enumerate}
    \item $\tr(\Op{T}_{u}\Op{T}_{u}^{*})=\tr(\Op{T}_{u}^{*}\Op{T}_{u})<\infty$.
    \item $u$ is cyclic, iff $\Op{T}_{u}$ is injective.
    \item $u$ is separating, iff $\Op{T}_{u}$ has dense range.
    \item $u$ is cyclic and separating iff $\Op{T}_{u}$ is injective and has
  dense range, i.e. iff $\Op{T}_{u}$ is invertible. 
  \end{enumerate}
\end{thm}
\begin{proof}
  The existence and the asserted properties follow from \lemref{lem2:1:1} and
  the uniqueness from \corref{cor2:1:1}. 
\end{proof}

\begin{thm}\label{thm2:2}
  Let $\Op{T}\eta\M_{0}$. Then $\Op{T}$ is the operator corresponding to a
  vector $u\in\H$ in the sense of \defnref{defn2:1:1} iff
  $\tr{\Op{T}\Op{T}^{*}}=\tr{\Op{T}^{*}\Op{T}}<\infty$.  
\end{thm}

\begin{proof}
  The necessarity of the trace condition follows from \thmref{thm2:1} and the
  sufficiency from \remref{rem2:1:1}.
\end{proof}

\section{Generation of Modular Objects}
\label{sec:3}

In this section we show how the modular objects of a cyclic and separating
vector $u_{0}\in\H$ for a semifinite von Neumann factor $(\M_{0},\H_{0})$ are
related to the operator $\Op{T}_{u_{0}}$ constructed in the last section.

Here again $\M_{0}$ is a finite factor with cyclic trace vector $u_{\tr}$
(cf.\secref{sec:2}). Then a conjugation $\Op{J}$ is defined by
\begin{equation}\label{eq3:1:1}
  \begin{split}
    \Op{J}:\H_{0}&\to\H_{0}\\
    \Op{A}u_{\tr}&\mapsto\Op{J}\Op{A}u_{\tr}=\Op{A}^{*}u_{\tr},
  \end{split}
\end{equation}
s.t. $\Op{A}\mapsto\Op{J}\Op{A}^{*}\Op{J}$ is an antiisomoprhism from $\M_{0}$
onto $\M_{0}^{'}$. Then we get the following
\begin{thm}\label{thm3:1}
  Let $\M_0$ be a finite von Neumann factor with cyclic and separating vector
  $u_{0}\in\M_{0}$ and cyclic trace vector $u_{\tr}\in\H_0$. Let further
  $\Op{T}_{u_{0}}\eta\M_{0}$ be the invertible operator corresponding to
  $u_{0}$ and
  $\Op{T}_{u_{0}}=\Op{HV}=(\Op{T}_{u_{0}}\Op{T}_{u_{0}}^{*})^{1/2}\Op{V}$ the
  polar decomposition of $\Op{T}_{u_{0}}$. Then we can calculate the modular
  objects $(\Delta_{0},\Op{J}_{0})$ of $(\M_{0},u_{0})$ as follows:
  \begin{equation*}
    \Op{J}_0=\Op{JV}^*\Op{JVJ}=\Op{VJV}^*,
  \end{equation*}
  where $\Op{J}$ is the
    conjugation defined in \eqref{eq3:1:1}, and
  \begin{equation*}
    \Delta_0=\Op{J}_0\Op{H}^{-1}_0\Op{J}_0\Op{H}_0,
  \end{equation*}
  where $\Op{H}_0=\Op{H}^2=\Op{T}_{u_{0}}\Op{T}_{u_{0}}^{*}$.
\end{thm}

\begin{proof}
  Let $\Op{E}_n:=\Op{E}_{\abge{1/n}{n}}\in\M_0$ be the spectral 
     projections of $\Op{H}$ corresponding to the interval ${\abge{1/n}{n}}$. 
     Then $\Op{T}_n:=(\Op{HE}_n+(\Op{I}-\Op{E}_n))\Op{V}$ is in 
     $\M_0$ with $\Op{V}$, $\Op{HE}_n+(\Op{I}-\Op{E}_n)$ as the polar 
     decomposition. According to \cite[9.6.11]{KRII} $u_n:=\Op{T}_n u_{\tr}$ is
     a sequence  of cyclic and separating vectors converging ($\Op{T}_n$ are
     invertible operators!) to $u$ with modular objects:
      \begin{equation*}
        \begin{split}
          \Op{J}_{n}&=\Op{JV}^*\Op{JVJ}=\Op{VJV}^*\text{ and }\\
          \Delta_{n}&=\Op{J}\Op{V}^{*}(\Op{HE}_n+(\Op{I}-\Op{E}_n))^{-2}
                        \Op{VJ}(\Op{HE}_n+(\Op{I}-\Op{E}_n))^2,
        \end{split}
      \end{equation*}
      where $\Op{J}$ is the conjugation corresponding to the trace vector
      $u_{\tr}$. 

      Now the modular conjugation of $u_{0}=\Lim{n}u_n$
      is $\Op{JV}^*\Op{JVJ}=\Op{VJV}^*$ since all $u_{n}$ lie in the same
      natural (closed) cone. Further the modular groups
      $\Delta_{n}^{it}$ corresponding to $u_n$ converge in the strong
      operator topology to the modular group $\Delta_0^{it}$ corresponding to
      $u_{0}$ (cf. \cite[p.106]{STR}). Since
      \begin{equation*}
        \Delta_{n}^{it}=
        \Op{J}\Op{V}^{*}(\Op{HE}_n+(\Op{I}-\Op{E}_n))^{2it}\Op{V}\Op{J}
        (\Op{HE}_n+(\Op{I}-\Op{E}_n))^{2it}
      \end{equation*}
      and operator multiplication is continuous on bounded sets w.r.t. the
      strong operator topology, we have 
      \begin{equation*}
        \begin{split}
          \Delta_u^{it}
          &=so-\Lim{n}\Delta_{n}^{it}\\
          &=\Op{J}\Op{V}^{*}\Op{H}^{2it}\Op{VJ}\Op{H}^{2it}\\
          &=\Op{VJV}^{*}\Op{H}^{2it}\Op{VJV}^{*}
            \Op{H}^{2it}\\
          &=\Op{J}_{u}\Op{H}_{0}^{it}\Op{J}_{u}\Op{H}_{0}^{it}.
        \end{split}
      \end{equation*}
      Since $\Op{J}_{u}\Op{H}_{0}^{-1}\Op{J}_{u}$ and $\Op{H}_{0}$ commute,
      $\Op{J}_{u}\Op{H}_{0}^{-1}\Op{J}_{u}\cdot\Op{H}_{0}$ is closable
      (cf. \cite[5.6.15]{KRI}) and the closure
      $\Op{J}_{u}\Op{H}_{0}^{-1}\Op{J}_{u}\Op{H}_{0}$ is selfadjoint,
      s.t. $\Delta_{u}=\Op{J}_{u}\Op{H}_{0}^{-1}\Op{J}_{u}\Op{H}_{0}$.
\end{proof}

\section{Two Simple Classes of Solutions of the Inverse Problem}
\label{sec:4}

In this section we want to use the results of the last two sections to examine
two simple classes of solutions of the inverse problem. How Wollenberg
showed in \cite{WolII} for arbitrary factors there is always a simple
class of solutions:
\begin{equation}\label{eq4:1}
  \begin{split}
    NF_{\M_{0}}^{1}(\Delta_{0},\Op{J}_{0},u_{0}):=\{&
    \M=\Op{U}\M_{0}\Op{U}^{*};\Op{U}\text{ unitary,}\\
    &\Op{U}u_{0}=\pm u_{0},\Op{UJ}_{0}=\Op{J}_{0}\Op{U},
    \Op{U}\Delta_{0}=\Delta_{0}\Op{U}\}.
  \end{split}
\end{equation}
Here $(\M_{0},\H_{0})$ is a von Neumann factor acting on a
Hilbert space $\H_{0}$ with a cyclic and separating vector $u_{0}$, and
$\Delta_{0}$, $\Op{J}_{0}$ are the modular objects w.r.t. $u_{0}$. For the
proof consider \propref{prop1:1}.

Further in \cite{WolII} was shown that, for finite type $I$ factors,
$(\Delta_{0}^{-1},\Op{J}_{0})$ are also modular objects for a cyclic and
separating vector $u_{1}$, and there is a unitary $\Op{U}_{1}$
s.t. $\Op{U}_{1}$ commutes with $\Op{J}_{0}$,
$\Delta_{0}^{-1}=\Op{U}_{1}^{*}\Delta_{0}\Op{U}_{1}$, and
$\Op{U}_{1}^{*}u_{j}=u_{j}$ ($j=0,1$). With this unitary a
second simple class of solutions was constructed:
\begin{equation}\label{eq4:2}
  \begin{split}
    NF_{\M_{0}}^{2}(\Delta_{0},\Op{J}_{0},u_{0}):=\{&
    \M=\Op{U}\M_{0}\Op{U}^{*};\Op{U}=\Op{KU}_{1},\Op{K}\text{ unitary,}\\
    &\Op{K}^{*}u_{0}=\pm u_{1},\Op{KJ}_{0}=\Op{J}_{0}\Op{K},
    \Op{K}\Delta_{0}=\Delta_{0}\Op{K}\}
  \end{split}
\end{equation}
Also here \propref{prop1:1} shows the assertion.

In the following we want to examine in the more general context of finite
factors, whether or not $\Delta_{0}^{-1}$ is also a modular operator for a
cyclic and separating vector, which then gives rise to a solution of the
inverse problem according to \propref{prop1:1}. For this purpose let
$\Delta_{0}=\Op{J}_{0}\Op{H}_{0}^{-1}\Op{J}_{0}\Op{H}_{0}$ be the decomposition
of the modular operator $\Delta_{0}$, where
$\Op{J}_0=\Op{JV}^*\Op{JVJ}=\Op{VJV}^*$ and
$\Op{T}_{u_{0}}=\Op{H}_{0}^{1/2}\Op{V}$ is the operator corresponding to
$u_{0}$ (cf. \thmref{thm3:1}). Then
$\Delta_{0}^{-1}=\Op{J}_{0}\Op{H}_{0}\Op{J}_{0}\Op{H}_{0}^{-1}$. For
$\Delta_{0}^{-1}$ being modular operator for a cyclic and separating vector it
is necessary and sufficient that $\tr(\Op{H}_{0}^{-1})<\infty$, which is shown
by the next

\begin{lem}\label{lem4:1}
  With the notations from above the following is equivalent:
  \begin{enumerate}
    \item ($\Delta_{0}^{-1},\Op{J}_{0})$ are the modular objects w.r.t. a
      cyclic and separating vector $u_{1}\in\H_{0}$.
    \item  
      \begin{equation}\label{eq4:3}
        \tr(\Op{H}_{0}^{-1})<\infty.
      \end{equation}
  \end{enumerate}
\end{lem}

\begin{proof}
  \begin{enumerate}
    \item "$\Rightarrow$":
      Suppose that $(\Delta_{0}^{-1},\Op{J}_{0})$ are the modular objects
      corresponding to a cyclic and separating vector $u_{1}$. Since $u_{1}$
      is cyclic and separating there exists a non-singular operator
      $\Op{S}\eta\M_0$ corresponding to $u_{1}$
      s.t. $\tr(\Op{S}\Op{S}^{*})=\tr(\Op{S}^{*}\Op{S})<\infty$ and
      \begin{equation*}
        \Delta_{0}^{-1}=\Op{SS}^{*}\Op{J}_{0}(\Op{SS}^{*})^{-1}\Op{J}_{0}
      \end{equation*}
      (cf. \thmref{thm2:1}, \thmref{thm2:2}, and \thmref{thm3:1}). Since the
      decomposition of $\Delta_{0}^{-1}$ is unique up to a positive constant
      (\propref{prop4:1}) we have 
      \begin{equation*}
        \tr(\Op{H}_{0}^{-1})=c\tr(\Op{S}\Op{S}^{*})<\infty.
      \end{equation*}
    \item "$\Leftarrow$":
      Suppose that $\tr(\Op{H}_{0}^{-1})<\infty$. Choosing the non-singular
      operator $\Op{S}:=\Op{H}_{0}^{-1/2}\Op{V}\eta\M_{0}$ (according to
      \remref{rem2:1:1} and \propref{prop2:1:2} $\Op{S}$ is densely defined
      and affiliated with $\M_{0}$) we have 
      \begin{equation*}
        \tr(\Op{S}\Op{S}^{*})=\tr(\Op{H}_{0}^{-1})<\infty
      \end{equation*}
      and
      \begin{equation*}
        \tr(\Op{S}^{*}\Op{S})=\tr(\Op{V}^{*}\Op{H}_{0}^{-1}\Op{V})
        =\tr(\Op{H}_{0}^{-1})<\infty.
      \end{equation*}
      Hence there is a cyclic and separating vector $u_{1}$ corresponding to
      $\Op{S}$ s.t. $(\Delta_{0}^{-1},\Op{J}_{0})$ are the modular objects
      w.r.t. $u_{1}$ (cf. \thmref{thm2:2} and \thmref{thm3:1}).
  \end{enumerate}
\end{proof}

Now we must examine, when the second condition in \lemref{lem4:1} is
fulfilled. For this purpose we must distinguish the types:
\begin{lem}\label{lem4:2}
  \begin{enumerate}
    \item For type $I_{N}$ factors ($N\in\N$) the condition \eqref{eq4:3} is
      always satisfied.
    \item For type $II_{1}$ factors the condition \eqref{eq4:3} is not always
      satisfied, but there are some vectors or operators,
      resp., s.t. the condition is fulfilled.
  \end{enumerate}
\end{lem}

\begin{proof}
  \begin{enumerate}
  \item 
    Let $\M_{0}$ be a type $I_{N}$ factor ($N\in\N$), i.e. it is
    isomorphic to $L(\H)$ with a finite dimensional $\H$. Since $\H$ is finite
    dimensional all linear operators are bounded and have finite trace,
    i.e. condition \eqref{eq4:3} is always satisfied.
  \item
    cf. \bspref{bsp4:1}.
 \end{enumerate}
\end{proof}

\begin{bsp}\label{bsp4:1}
  Let $\H:=L_{2}(S,\alg{S},m)$ and $\alg{A}:=L_{\infty}(S,\alg{S},m)$ be the
  multiplication algebra on $\H$, where $S=\offenr{0}{1}$, $\alg{S}$ the Borel
  sets in $S$, and $m$ the Lebesgue measure. Let further $G$ be the group of
  all rational translations, modulo 1, of $S$. Then the crossed product von
  Neumann algebra $\alg{R}$ generated by $\alg{A}$ and $G$ is the hyperfinite
  factor of type $II_{1}$ (cf. \cite[Ex.8.6.12]{KRII}). Let
  $\Op{M}=(\Op{U}(pq^{-1})\Op{A}(pq^{-1}))_{p,q\in G}\in\alg{R}$, where
  $\Op{U}$ is the unitary representation of $G$ on $\H$, and
  $\Op{A}(p)\in\alg{A}$ for all $p\in G$, then 
  \begin{equation*}
    \tr{\Op{M}}:=\int f dm
  \end{equation*}
  with $f\in L_{\infty}$ s.t. $\Op{M}_{f}=\Op{A}(e)$, is the trace on
  $\alg{R}$. Let $\Phi$ be the canonical isomorphism from $\alg{A}$ into
  $\alg{R}$ ($\Phi(f)=(\delta_{pq}\Op{M}_{f})_{p,q}$), then
  \begin{equation*}
    \tr(\Phi(f)^{*}\Phi(f))<\infty\quad\Leftrightarrow\quad
    \int\abs{f}^{2}dm<\infty.
  \end{equation*}
  Now let $f_{1}:=x+1$ and $f_{2}:=x$. Then $f_{1},f_{2}\in\alg{A}$ and
  \begin{equation*}
    \int_{0}^{1}(x+1)^{2}dm(x)<\infty\text{ and }
    \int_{0}^{1}(x+1)^{-2}dm(x)<\infty
   \end{equation*}
   whereas
  \begin{equation*}
    \int_{0}^{1}x^{2}dm(x)<\infty\text{ and }
    \int_{0}^{1}x^{-2}dm(x)=\infty.
   \end{equation*}
   This shows that for $\Phi(f_{1})$ condition \eqref{eq4:3} is satisfied
   whereas for $\Phi(f_{2})$ not.

   Since every type $II_{1}$ factor has a hyperfinite subfactor
   (cf. \cite[Ex.12.4.25]{KRII}) this example
   also shows the second assertion of \lemref{lem4:2} for all type $II_{1}$
   factors. 
\end{bsp} 

\begin{rem}\label{rem4:0}
   In \bspref{bsp4:1} we can observe that the spectral measure of $\Phi(f)$ is
   \begin{equation*}
     \Op{E}_{\Phi(f)}(B)=(\Op{M}_{f^{-1}(B)}\delta_{pq})_{p,q}.
   \end{equation*}
   This shows that $\Phi(f)$ has the same spectrum as $f$, i.e. all types of
   spectral points can appear for the positive operator $\Op{H}_{0}$ generating
   the modular operator - in contrast to the type $I$ case, where we
   have, also for $I_{\infty}$, only point spectrum (and $0$ as continuous
   spectrum), since there $\Op{H}_{0}$ is always a trace class operator.
\end{rem}

Suppose in the following that one (and hence all) of the conditions of
\lemref{lem4:1} is true. According to
\remref{rem2:1:1} then for $\Op{S}:=\Op{H}_{0}^{-1/2}\Op{V}$ with
$\tr(\Op{S}^{*}\Op{S})=\tr(\Op{H}_{0}^{-1})<\infty$ (for the definition of
$\Op{S}$ see the proof of \lemref{lem4:1}) $u_{\tr}\in\Dom{S}$, and
$u_{1}:=\Op{H}_{0}^{-1/2}\Op{V}u_{\tr}$ is a cyclic and separating vector with
$(\Delta_{0}^{-1},\Op{J}_{0})$ as modular objects. Notice further that in this
case $\tr(\Op{H}_{0}^{-1/2})<\infty$ also is equivalent to
$u_{\tr}\in\Dom{H_{0}^{-1/2}}$. Then we can
construct with the help of
\lemref{lem3:1} a conjugation $\Op{I}$, s.t. $\Op{I}$ commutes with
$\Delta_{0}$ and $\Op{J}_{0}$, and $\Op{I}u_{i}=u_{i}$ ($i=0,1$). Setting
$\Op{U}_{1}:=\Op{I}\Op{J}_{0}$ we get unitary commuting with $\Op{J}_{0}$,
$\Op{U}_{1}u_{i}=u_{1}$, and
\begin{equation*}
  \Op{U}_{1}^{*}\Delta_{0}\Op{U}_{1}=\Delta_{0}^{-1}.
\end{equation*}
Now we can define the following class of von Neumann factors solving the
inverse problem (cf. \eqref{eq4:2}):
\begin{equation}\label{eq4:4}
  \begin{split}
    NF_{\M_0}^{2}(\Delta_{0},\Op{J}_{0},u_{0}):=\{&
    \M=\Op{U}\M_{0}\Op{U}^{*};\Op{U}=\Op{KU}_{1},\Op{K}\text{ unitary,}\\
    &\Op{K}^{*}u_{0}=\pm u_{1},\Op{KJ}_{0}=\Op{J}_{0}\Op{K},
    \Op{K}\Delta_{0}=\Delta_{0}\Op{K}\}
  \end{split}
\end{equation}
and state the following
\begin{lem}\label{lem4:3}
  Let $\M_{0}$ be a finite von Neumann factor,
  $u_{0}=\Op{H}_{0}^{-1/2}\Op{V}u_{\tr}$ a
  cyclic and separating vector for $\M_{0}$,
  s.t. $u_{\tr}\in\Dom{H_{0}^{-1/2}}$, or, equivalently,
  $u_{\tr}\in\Dom{H_{0}^{-1/2}V}$. Then
  \begin{equation*}
    NF_{\M_0}^{2}(\Delta_{0},\Op{J}_{0},u_{0})\subset
    NF_{\M_0}(\Delta_{0},\Op{J}_{0},u_{0})
  \end{equation*}
\end{lem}

The proof is the same as for type I factors (s. \cite{WolII}). Also in this case
we can show the next lemma and its corollary analogous to Prop. 3.4 and
Prop. 3.5 in \cite{WolII}.

\begin{lem}\label{lem4:4}
  Let $\M=\Op{U}\M_{0}\Op{U}^{*}\in
  NF_{\M_0}(\Delta_{0},\Op{J}_{0},u_{0})$, where $\Op{U}$ commutes with
  $\Op{J}_{0}$. Let further
  $(\Delta:=\Op{U}^{*}\Delta_{0}\Op{U},\Op{J}_{0})$ be the modular objects for
  $(\M_{0},\Op{U}^{*}u_{0})$ and
  $\Delta_{0}=\Op{H}_{0}\Op{J}_{0}\Op{H}_{0}^{-1}\Op{J}_{0}$,
  $\Delta=\Op{H}\Op{J}_{0}\Op{H}^{-1}\Op{J}_{0}$ with
  $\Op{H}_{0},\Op{H}\eta\M_{0}$. Then
  \begin{enumerate}
    \item $\M\in NF_{\M_0}^{1}$ iff there exists a unitary
      $\Op{W}_{1}\in L(\H)$, s.t. $\Op{H}_{0}=c\Op{W}_{1}\Op{HW}_{1}^{*}$ and
      $\ad\Op{W}_{1}\in \aut\M_{0}$ and $\Op{W}_{1}$ commutes with
      $\Op{J}_{0}$. 
    \item $\M\in NF_{\M_0}^{2}$ iff there exists a unitary
      $\Op{W}_{2}\in L(\H)$,
      s.t. $\Op{H}_{0}=c\Op{W}_{2}\Op{H}^{-1}\Op{W}_{2}^{*}$ and
      $\ad\Op{W}_{2}\in \aut\M_{0}$ and $\Op{W}_{2}$ commutes with
      $\Op{J}_{0}$. 
  \end{enumerate}
\end{lem}

\begin{proof}
  \begin{enumerate}
    \item "$\Rightarrow$": Let $\M=\Op{U}\M_{0}\Op{U}^{*}\in
      NF_{\M_0}^{1}$. Then there exists a unitary $\Op{V}\in L(\H)$ which
      commutes with $\Delta_{0}$ and $\Op{J}_{0}$,
      s.t. $\M=\Op{V}\M_{0}\Op{V}^{*}$ and $\Op{V}^{*}u_{0}=\pm u_{0}$. Setting
      $\Op{W}_{1}:=\Op{V}^{*}\Op{U}$ we have $\ad \Op{W}_{1}\in\aut\M_{0}$ and
      $\Op{W}_{1}$ commutes with $\Op{J}_{0}$ (since $\Op{U}$ and $\Op{V}$
      do). Then we can calculate
      \begin{equation}\label{eq4:5}
        \begin{split}
          \Op{W}_{1}\Delta\Op{W}_{1}^{*}
          &=\Op{V}^{*}\Op{U}\Op{U}^{*}\Delta_{0}\Op{UU}^{*}\Op{V}\\
          &=\Op{V}^{*}\Delta_{0}\Op{V}\\
          &=\Delta_{0}=\Op{H}_{0}\Op{J}_{0}\Op{H}_{0}^{-1}\Op{J}_{0}
        \end{split}
      \end{equation}
      Since $\Op{W}_{1}$ commutes with $\Op{J}_{0}$ we have
      \begin{equation*}
        \Op{W}_{1}\Delta\Op{W}_{1}^{*}=
        (\Op{W}_{1}\Op{H}\Op{W}_{1}^{*})
        (\Op{J}_{0}\Op{W}_{1}\Op{H}^{-1}\Op{W}_{1}^{*}\Op{J}_{0}).
      \end{equation*}
      With this, \eqref{eq4:5}, and \lemref{lem4:1} follows
      \begin{equation*}
        \Op{H}_{0}=c\Op{W}_{1}\Op{H}\Op{W}_{1}^{*}.
      \end{equation*}

      "$\Leftarrow$": We can assume w.l.o.g. that
      \begin{equation*}
        \Op{H}_{0}=\Op{W}_{1}\Op{H}\Op{W}_{1}^{*},
      \end{equation*}
      where $\Op{W}_{1}\in L(\H_{0})$ unitary and $\ad
      \Op{W}_{1}\in\aut\M_{0}$ and $\Op{W}_{1}$ commutes with
      $\Op{J}_{0}$. Then
      \begin{equation*}
        \begin{split}
          \Op{W}_{1}^{*}\Delta_{0}\Op{W}_{1}&=
          \Op{W}_{1}^{*}\Op{H}_{0}\Op{W}_{1}
          \Op{J}_{0}\Op{W}_{1}^{*}\Op{H}_{0}^{-1}\Op{W}_{1}\Op{J}_{0}\\
          &=\Op{H}\Op{J}_{0}\Op{H}^{-1}\Op{J}_{0}=\Delta.
        \end{split}
      \end{equation*}
      Since $\Delta=\Op{U}^{*}\Delta_{0}\Op{U}$, we define
      $\Op{K}:=\Op{U}\Op{W}_{1}^{*}$, s.t. $\Op{K}$ commutes with $\Delta_{0}$
      and 
      $\Op{J}_{0}$. Further $u_{0}$ is cyclic and separating for
      $\M=\Op{K}\M_{0}\Op{K}^{*}=\Op{U}\M_{0}\Op{U}^{*}$ and
      $(\Delta_{0},\Op{J}_{0})$ are the modular objects for
      $(\M_{0},u_{0})$. This means that
      $(\Op{K}^{*}\Delta_{0}\Op{K}=\Delta_{0},\Op{K}^{*}\Op{J}_{0}\Op{K}=\Op{J}_{0})$ are the modular objects for $(\M_{0},\Op{K}^{*}u_{0})$. But $(\Delta_{0},\Op{J}_{0})$ are the modular objects for $(\M_{0},u_{0})$, too. From this follows, that $\Op{K}^{*}u_{0}=\pm u_{0}$, and $\M\in NF_{\M_0}^{1}$.
    \item Analogous to part 1, cf. \cite[Prop. 3.5]{WolII}.
  \end{enumerate}
\end{proof}

\begin{cor}\label{cor4:1}
  Let $\Delta_{0}=\Op{H}_{0}\Op{J}_{0}\Op{H}_{0}^{-1}\Op{J}_{0}$ the
  canonical decomposition of the modular operator $\Delta_{0}$. Suppose 
  \begin{equation*}
    \M\in NF_{\M_0}^{1}(\Delta_{0},\Op{J}_{0},u_{0})\cap 
          NF_{\M_0}^{2}(\Delta_{0},\Op{J}_{0},u_{0}).
  \end{equation*}
  Then $\Op{H}_{0}=c\Op{V}\Op{H}_{0}^{-1}\Op{V}^{*}$, where $\Op{V}$ is a
  unitary in $\H_{0}$, $\ad\Op{V}\in\aut\M_{0}$, $c>0$ and $\Op{V}$ commutes
  with $\Op{J}_{0}$.
\end{cor}

\begin{proof}
  Let $\Op{U}$ be a unitary, s.t. $\M=\Op{U}\M_{0}\Op{U}^{*}$ and 
  \begin{equation*}
    \Delta:=\Op{U}^{*}\Delta_{0}\Op{U}=\Op{H}\Op{J}_{0}\Op{H}^{-1}\Op{J}_{0},
  \end{equation*}
  where $\Op{H}\eta\M_{0}$. Then from \lemref{lem4:4} follows
  \begin{gather*}
    \Op{H}_{0}=c_{1}\Op{W}_{1}\Op{H}\Op{W}_{1}^{*}\text{ and}\\
    \Op{H}_{0}=c_{2}\Op{W}_{2}\Op{H}^{-1}\Op{W}_{2}^{*}.
  \end{gather*}
  This shows 
  \begin{equation*}
    \Op{H}_{0}=
    c_{1}c_{2}^{-1}\Op{W}_{1}\Op{W}_{2}^{*}\Op{H}_{0}^{-1}\Op{W}_{2}\Op{W}_{1}^{*}. 
  \end{equation*}
\end{proof}

\begin{rem}\label{rem4:1}
  In the case of type $I$ factors Wollenberg showed in \cite{WolII} that there
  is a special class of modular operators with so called generic spectrum
  s.t.  $NF_{\M_0}^{1}$ and $NF_{\M_0}^{2}$ are the only classes of solutions
  of the inverse problem. As it will be shown in the next chapter, for modular
  operators generated by operators with pure point spectrum not being the
  identity operator, in the type $II$ case there are always more classes. 
\end{rem}

\section{The Classification of Solutions in the Pure Point Spectrum Case}
\label{sec:5}
If the non-singular, positive operator $\Op{H}_{0}$ which generates the
modular operator
$\Delta_{0}=\Op{H}_{0}\Op{J}_{0}\Op{H}_{0}^{-1}\Op{J}_{0}$ has pure point
spectrum (and then $\Delta_{0}$ have it, too), we can give a complete
classification of solutions of the inverse problem in terms of the spectrum of
$\Op{H}_{0}$. For this purpose we define the following equivalence relation
in the set of solutions of the inverse problem:
\begin{defn}\label{defn5:1}
  Two semifinite von Neumann factors
  $\M,\alg{N}\in NF_{\M_{0}}(\Delta_{0},\Op{J}_{0},u_{0})$ are called
  equivalent, $\M\sim\alg{N}$, if $\M\in
  NF^{1}_{\alg{N}}(\Delta_{0},\Op{J}_{0},u_{0})$, i.e. if there exists
  a unitary operator $\Op{U}$ on $\H_{0}$,
  s.t. $\M=\Op{U}\alg{N}\Op{U}^{*}$, $\Op{U}$ commutes with $\Delta_{0}$ and
  $\Op{J}_{0}$ and $\Op{U}^{*}u_{0}=\pm u_{0}$ (For the definition of the
  class $NF^{1}_{\alg{N}}(\Delta_{0},\Op{J}_{0},u_{0})$ see \eqref{eq4:1}).
  \end{defn}

  \begin{proof}
    The relation defined in \defnref{defn5:1} is an equivalence relation, since
    it is reflexive (choose $\Op{U}=\Op{Id}$), symmetric
    ($\tilde\Op{U}=\Op{U}^{*}$) and transitive: Let $\M\sim\alg{N}$ and
    $\alg{N}\sim\alg{R}$; this means hat $\M=\Op{U}\alg{N}\Op{U}^{*}$ and
    $\alg{N}=\Op{V}\alg{R}\Op{V}^{*}$, where $\Op{U}$ and $\Op{V}$ are
    unitaries with the properties described above. Then with
    $\tilde\Op{U}:=\Op{UV}$ we have $\M=\tilde\Op{U}\alg{R}\tilde\Op{U}^{*}$
    and also $\tilde\Op{U}$ have the right properties, s.t. $\M\sim\alg{R}$. 
  \end{proof}

\begin{rem}\label{rem5:1}
  \begin{enumerate}
    \item It is trivial from the definition of the equivalence relation that
      the first simple class
      $NF^{1}_{\alg{M}_{0}}(\Delta_{0},\Op{J}_{0},u_{0})$ of solutions of the
      inverse problem is an equivalence class w.r.t. this equivalence
      relation.
    \item Also the second simple class from \secref{sec:4} is an equivalence
      class w.r.t. $\sim$.
      \begin{proof}
        Let $M_{i}\in NF_{\M_{0}}^{2}(\Delta_{0},\Op{J}_{0},u_{0})$ $(i=1,2)$
        be two members of this class. Then there exist unitaries $\Op{K}_{i}$,
        s.t. $\Op{K}_{i}^{*}u_{0}=\pm u_{1}$, $\Op{K}_{i}$ commutes with
        $\Op{J}_{0}$ and $\Delta_{0}$, and
        $\M_{i}=\Op{K}_{i}\Op{U}_{1}\M_{0}\Op{U}_{1}^{*}\Op{K}_{i}^{*}$
        (cf. \eqref{eq4:4}). Now define
        \begin{equation*}
          \Op{W}:=\Op{K}_{1}\Op{U}_{1}\Op{U}_{1}^{*}\Op{K}_{2}^{*}
          =\Op{K}_{1}\Op{K}_{2}^{*}.
        \end{equation*}
        Then $\M_{1}=\Op{W}\M_{2}\Op{W}^{*}$ and
        \begin{gather*}
          \Op{W}\Delta_{0}=\Op{K}_{1}\Op{K}_{2}^{*}\Delta_{0}
          =\Delta_{0}\Op{K}_{1}\Op{K}_{2}^{*}=\Delta_{0}\Op{W}\\
          \Op{W}\Op{J}_{0}=\Op{K}_{1}\Op{K}_{2}^{*}\Op{J}_{0}
          =\Op{J}_{0}\Op{K}_{1}\Op{K}_{2}^{*}=\Op{J}_{0}\Op{W}\\
          \Op{W}^{*}u_{0}=\Op{K}_{2}\Op{K}_{1}^{*}u_{0}
          =\pm\Op{K}_{2}u_{1}=\pm u_{0},
        \end{gather*}
        s.t. the conditions of \defnref{defn5:1} are fulfilled and
        $\M_{1}\sim\M_{2}$.
      \end{proof}
  \end{enumerate}
\end{rem}

Let now $(\Delta_{0},\Op{J}_{0})$ be the modular objects of $(\M_{0},u_{0}$),
where $\Op{T}_{u_{0}}=\Op{H}_{0}^{1/2}\Op{V}$ is the non-singular operator
corresponding to the cyclic and separating vector $u_{0}\in\H_{0}$ and 
$\Delta_{0}:=\Op{H}_{0}\Op{J}_{0}\Op{H}_{0}\Op{J}_{0}$ ($(\M_{0},\H_{0})$ is a
finite von Neumann factor). If a factor $\M$ is a solution of the inverse
Problem, $\M\in NF_{\M_{0}}(\Delta_{0},\Op{J}_{0},u_{0})$, then there is,
according to \propref{prop1:1}, a unitary $\Op{U}$,
s.t. $(\Delta:=\Op{U}^{*}\Delta_{0}\Op{U},\Op{J}_{0})$ are the modular objects
for $(\M_{0},u:=\Op{U}^{*}u_{0})$. According to \thmref{thm3:1} there is a
positive non-singular operator $\Op{H}\eta\M_{0}$,
s.t. $\Delta=\Op{H}\Op{J}_{0}\Op{H}^{-1}\Op{J}_{0}$. \thmref{thm2:1} shows
$\tr(\Op{H})<\infty$, whence we can assume w.l.o.g. 
\begin{equation}\label{eq5:1}
  \tr(\Op{H})=1.
\end{equation}
What does the equivalence relation defined in \defnref{defn5:1} mean for these
operators $\Op{H}$? This question is answered by the next
\begin{lem}\label{lem5:0}
  Let $\M_{1},\M_{2}\in NF_{\M_{0}}(\Delta_{0},\Op{J}_{0},u_{0})$ be two
  solutions of the inverse problem, and $\Op{H}_{i}\eta\M_{0}$ ($i=1,2$) the
  corresponding positive operators. Then the following is equivalent:
  \begin{enumerate}
  \item $\M_{1}\sim\M_{2}$
  \item There is a unitarily implemented automorphism $\alpha=\ad W$,
    $\Op{W}\in\alg{U}(\H_{0})$, of
    $\M_{0}$, s.t. $\alpha(\Op{H}_{1})=\Op{H}_{2}$,
    i.e. $\Op{W}\Op{H}_{1}\Op{W}^{*}=\Op{H}_{2}$.
  \end{enumerate}
  The operator $\Op{W}$ can be chosen in such a way that it commutes with
  $\Op{J}_{0}$.
\end{lem}

\begin{proof}
  Since $\M_{1}$ and $\M_{2}$ are solutions of the inverse problem, there are
    two unitary operators $\Op{U}_{1}$ and $\Op{U}_{2}$ both commuting with
    $\Op{J}_{0}$, s.t. $\M_{i}=\Op{U}_{i}\M_{0}\Op{U}_{i}^{*}$,
    $(\Delta_{i}:=\Op{U}_{i}^{*}\Delta_{0}\Op{U}_{i}
    =:\Op{H}_{i}\Op{J}_{0}\Op{H}_{i}^{-1}\Op{J}_{0},\Op{J}_{0})$ are the
    modular objects for $(\M_{0},u_{i}:=\Op{U}_{i}^{*}u_{0})$ ($i=1,2$).
  \begin{enumerate}
  \item Let $\M_{1}\sim\M_{2}$. Then there is a
    unitary $\Op{V}$ commuting with $\Op{J}_{0}$ and $\Delta_{0}$,
    s.t. $\M_{1}=\Op{V}\M_{2}\Op{V}^{*}$. Setting
    $\Op{W}:=\Op{U}_{1}^{*}\Op{V}\Op{U}_{2}$ an easy calculation gives 
    $\ad W\in\aut(\M_{0})$ and $\Op{W}$ commutes with $\Op{J}_{0}$. Also we
    can calculate 
    \begin{equation*}
      \begin{split}
        (\Op{W}\Op{H}_{2}\Op{W}^{*})
        (\Op{J}_{0}\Op{W}\Op{H}_{2}^{-1}\Op{W}^{*}\Op{J}_{0})
        &=\Op{U}_{1}^{*}\Op{V}\Delta_{0}\Op{V}\Op{U}_{1}\\
        &=\Op{H}_{1}\Op{J}_{0}\Op{H}_{1}^{-1}\Op{J}_{0}.
      \end{split}
    \end{equation*}
    \propref{prop4:1} now shows that $\Op{W}\Op{H}_{2}\Op{W}^{*}=\Op{H}_{1}$
    (note our normalization condition \eqref{eq5:1}).
  \item Suppose now $\alpha\in\aut(\M_{0})$. Since $\M_{0}$ possesses a cyclic
    and separating vector there is a unitary $\Op{W}$,
    s.t. $\alpha=\ad\Op{W}$. Further
    $(\Op{W}\M_{0}\Op{W}^{*}=\M_{0},\Op{W}u_{1})$ has modular objects
    $(\Op{W}\Delta_{1}\Op{W}^{*},\Op{W}\Op{J}_{0}\Op{W}^{*})$. Now there is a
    cyclic and separating vector $v$ in the natural cone of $u_{0}$ and a
    unitary $\Op{U}^{'}\in\alg{U}(\M_{0}^{'})$,
    s.t. $\Op{W}u_{1}=\Op{U}^{'}v$, and
    $(\M_{0},v)=(\M_{0},{\Op{U}^{'}}^{*}\Op{W}u_{1})$ has modular objects 
    $(\Op{W}\Delta_{1}\Op{W}^{*},\Op{J}_{0})=
     (\Op{W}\Delta_{1}\Op{W}^{*},
      {\Op{U}^{'}}^{*}\Op{W}\Op{J}_{0}\Op{W}^{*}\Op{U}^{'})$,
    i.e. $\tilde\Op{W}:={\Op{U}^{'}}^{*}\Op{W}$ commutes with $\Op{J}_{0}$
    and, since $\Op{U}^{'}\in\M_{0}^{'}$, $\ad\Op{W}=\ad\tilde\Op{W}$. Now
    define $\Op{U}:=\Op{U}_{1}{\tilde\Op{W}}^{*}\Op{U}_{2}^{*}$. Then
    $\Op{U}$ commutes with $\Op{J}_{0}$ and $\Delta_{0}$, for:
    \begin{equation*}
      \begin{split}
        \Op{U}^{*}\Delta_{0}\Op{U}
        &=\Op{U}_{2}\tilde\Op{W}\Op{U}_{1}^{*}\Delta_{0}
          \Op{U}_{1}{\tilde\Op{W}}^{*}\Op{U}_{2}^{*}\\
        &=\Op{U}_{2}\tilde\Op{W}\Op{H}_{1}\Op{J}_{0}\Op{H}_{1}^{-1}
          \Op{J}_{0}{\tilde\Op{W}}^{*}\Op{U}_{2}^{*}\\
        &=\Op{U}_{2}\Op{H}_{2}\Op{J}_{0}\Op{H}_{2}^{-1}
          \Op{J}_{0}\Op{U}_{2}^{*}\\
        &=\Delta_{0}.
      \end{split}
    \end{equation*}
    Also $\Op{U}\M_{2}\Op{U}^{*}=\M_{1}$, and, since $(\M_{0},u_{1})$ has
    modular objects $(\Delta_{1},\Op{J}_{0})$,
    $(\tilde\Op{W}\M_{0}{\tilde\Op{W}}^{*}=\M_{0},
      \tilde\Op{W}u_{1})$ has modular objects
    $(\tilde\Op{W}\Delta_{0}{\tilde\Op{W}}^{*}\Op{U}^{'}=\Delta_{2},
      \Op{J}_{0})$. Now, since the cyclic and separating vector is (up to the
    sign) uniquely determined by the modular objects
    (s. \cite[ch.2.(i)]{WolII}), $\tilde\Op{W}u_{1}=\pm u_{2}$,
    i.e. $\Op{U}^{*}u_{0}=\pm\Op{U}_{2}u_{2}=\pm u_{0}$, and $\Op{U}$ is the
    unitary required by \defnref{defn5:1}.
  \end{enumerate}
\end{proof} 

The last lemma says that for classifying the equivalence classes of $\sim$ we
must search for a complete set of invariants of selfadjoint operators  under
automorphisms. In general such a set is not known. But if we have an operator
with pure point spectrum, we can give such a set. Thus assume in
the following that $\Op{H}_{0}$ has pure point spectrum, i.e. 
$\Op{H}_{0}=\sum_{k\in K}\mu_{k}\Op{E}_{k}$ where the $\mu_{k}$ ($k\in K$) are
the eigenvalues of
$\Op{H}_{0}$ and $\Op{E}_{k}\in\M_{0}$ are the corresponding (orthogonal)
eigenprojections with
$m_{k}:=\tr{\Op{E}_{k}}=:D_{\M_{0}}(\Op{E}_{k})$ their von Neumann
dimension (cf. \cite[8.4]{KRII} for the notion of dimension in von Neumann
factors).

Now we have for $\Delta_{0}$ the following decomposition
\begin{equation}\label{eq5:3}
  \begin{split}
    \Delta_{0}&=\Op{H}_{0}\Op{J}_{0}\Op{H}_{0}^{-1}\Op{J}_{0}\\
              &=\sum_{k,l\in K}\mu_{k}\mu_{l}^{-1}\Op{E}_{k}
                         \Op{J}_{0}\Op{E}_{l}\Op{J}_{0}\\
              &=\sum_{j\in J}\lambda_{j}\Op{F}_{j},
  \end{split}
\end{equation}
where the $\lambda_{j}$ ($j\in J$) are the eigenvalues of $\Delta_{0}$, which
are invariant under unitary transformations $\Op{U}\in L(\H_{0})$, and
$\Op{F}_{j}$ are the corresponding eigenprojections. Now we can formulate the
following

\begin{lem}\label{lem5:1} 
  With the notations introduced above we can compute the spectrum of
  $\Delta_{0}$  in the following way: 
  \begin{equation}\label{eq5:4}
    \{\lambda_{j}\vert j\in J\}=\{\mu_{k}\mu_{l}^{-1}\vert k,l\in K\}
    \quad\forall j\in J
  \end{equation}
  and
  \begin{subequations}\label{eq5:5}
    \begin{equation}\label{eq5:5a}
      n_{j}=\sum_{\mu_{k}\mu_{l}^{-1}=\lambda_{j}}m_{k}m_{l}
      \quad\forall j\in J\text{ if $\M_{0}$ is type $I$,}
    \end{equation}
    \begin{equation}\label{eq5:5b}
      n_{j}=\infty\quad\forall j\in J\text{ if $\M_{0}$ is type $II$,}
    \end{equation}  
  \end{subequations}
  where $n_{j}:=D_{L(\H_{0})}(\Op{F}_{j})$ with $D_{L(\H_{0})}(\Op{F}_{j})$
  the dimension function in
  the type $I_{\infty}$ factor $L(\H_{0})$, which corresponds to the
  normalized Hilbert space dimension.
\end{lem}

For the proof we need the following
\begin{prop}\label{prop5:5}
  Let $\Op{E}\in\alg{P}(\M)$ and $\Op{F}\in\alg{P}(\M^{'})$ be two projections,
  where $\alg{P}(\M)$ and $\alg{P}(\M^{'})$ are the sets of projections in a
  non type $I$ von Neumann factor ($\M,\H)$ and its commutant, resp. Then the
  product $\Op{E}\Op{F}\in L(H)$ has infinite Hilbert space dimension,
  i.e. $D(\Op{EF})=\infty$.
\end{prop}

\begin{proof}
  Since $\M$ and $\M^{'}$ has no non-zero Abelian projection there exist for
  every $N\in\N$ orthogonal families 
  $(\Op{E}_{n}^{(N)})_{1\leq n\leq N}\subset\M$ and 
  $(\Op{F}_{n}^{(N)})_{1\leq n\leq N}\subset\M^{'}$ of non-zero projections,
  s.t. 
  \begin{gather*}
    \sum_{n=1}^{N}\Op{E}_{n}^{(N)}=\Op{E}\\
    \text{and}\\
    \sum_{n=1}^{N}\Op{F}_{n}^{(N)}=\Op{F}
  \end{gather*}
  (cf. \cite[Lemma 6.5.6]{KRII}). Now
  \begin{equation*}
    \Op{E}\Op{F}=\sum_{n,m=1}^{N}\Op{E}_{n}^{(N)}\Op{F}_{m}^{(N)},
  \end{equation*}
  where the projections $\Op{E}_{n}^{(N)}\Op{F}_{m}^{(N)}$ are pairwise
  orthogonal projections, since $\Op{E}_{n}^{(N)}$ commutes with
  $\Op{F}_{m}^{(N)}$, they are not $0$, since
  $\M$ is a factor (cf. \cite[prop.5.5.3]{KRI}), the latter means that they
  have at least Hilbert space dimension
  $1$. This means that the dimension of $\Op{E}\Op{F}$ is at least $N^{2}$ and,
  since $N$ was arbitrary, infinite. 
\end{proof}

\begin{proof}[Proof of \lemref{lem5:1}]
  The first assertion
  follows directly from \eqref{eq5:3} and the fact, that
  \begin{equation*}
    \Op{F}_{j}=
               \sum_{\mu_{k}\mu_{l}^{-1}=\lambda_{j}}
                   \Op{E}_{k}\Op{J}_{0}\Op{E}_{l}\Op{J}_{0}\not=0,
  \end{equation*}
  since $\M_{0}$ is a factor (cf. \cite[prop.5.5.3]{KRI}). The second assertion
  follows from $L(\H_{0})=\M_{0}\otimes\M_{0}^{'}$ and
  $\tr_{L(\H_{0})}=\tr_{\M_{0}}\otimes\tr_{M_{0}^{'}}$ in the type $I$ case
  and from \propref{prop5:5} in the type $II$ case.
\end{proof}

The next proposition shows that the eigenvalues and multiplicities characterize
a given operator with pure point spectrum in a von Neumann factor uniquely up
to unitary equivalence in the von Neumann factor.

\begin{prop}\label{prop5:2}
  Let $\Op{H},\tilde\Op{H}\eta\M_{0}$ be two selfadjoint operators with pure
  point spectrum affiliated
  with a
  semifinite von Neumann factor $\M_{0}$ which have the same eigenvalues and
  von Neumann multiplicities w.r.t. $\M_{0}$. Then there is a unitary
  $\Op{W}\in\M_{0}$ s.t. $\tilde\Op{H}=\Op{W}\Op{H}\Op{W}^{*}$.  
\end{prop}

\begin{proof}
    Since $\Op{H},\tilde\Op{H}\eta\M_{0}$ are two selfadjoint operators having
    the same eigenvalues $\mu_{k}$ ($k\in K$) we can write
    \begin{gather*}
      \Op{H}=\sum_{k\in K}\mu_{k}\Op{E}_{k}\\
      \tilde\Op{H}=\sum_{k\in K}\mu_{k}\Op{F}_{k}
    \end{gather*}
    where $\Op{E}_{k},\Op{F}_{k}\in\M_{0}$ are the corresponding (orthogonal)
    eigenprojections (and $\sum_{k\in K}\Op{E}_{k}=\sum_{k\in
      K}\Op{F}_{k}=\Op{Id}$) and the convergence is understood in the
    so-topology. Since $\Op{H}$ and $\tilde\Op{H}$ have the same
    multiplicities we have $D_{\M_{0}}(\Op{E}_{k})=D_{\M_{0}}(\Op{F}_{k})$
    ($k\in K$), where $D$
    is the unique dimension function on $\M_{0}$. This means that there are
    partial isometries $\Op{W}_{k}$, s.t. $\Op{W}_{k}^{*}\Op{W}_{k}=\Op{E}_{k}$
    and $\Op{W}_{k}\Op{W}_{k}^{*}=\Op{F}_{k}$
    (cf. \cite[Th.8.4.3]{KRII}). Setting $\Op{W}:=so-\sum_{k\in K}\Op{W}_{k}$
    we get a unitary in $\M_{0}$, s.t. $\Op{W}\Op{E}_{k}\Op{W}^{*}=\Op{F}_{k}$
    ($k\in K$) and $\Op{W}\Op{H}\Op{W}^{*}=\tilde\Op{H}$.
\end{proof}

Now we can show the following 
\begin{lem}\label{lem5:2}
  If there are two solutions of the inverse problem $\M_{1}$, $\M_{2}$ s.t. the
  corresponding selfadjoint operators $\Op{H}_{1}$ and $\Op{H}_{2}$ have the
  same eigenvalues modulo a positive constant $c>0$ and same (von Neumann)
  multiplicities, then $\M_{1}\sim\M_{2}$. 
\end{lem}

\begin{proof}
  According to \lemref{lem5:0} we have to show that there is an automorphism
  $\alpha=\ad\Op{W}\in\aut(\M_{0})$
  s.t. $\Op{H}_{2}=\Op{W}\Op{H}_{1}\Op{W}^{*}$. But for $\Op{W}$ we can chose
  the operator existing according to \propref{prop5:2}.
\end{proof}

The converse shows the next
\begin{lem}\label{lem5:3}
  If there are two equivalent solutions $\M_{1}$, $\M_{2}$ of the inverse
  problem with the corresponding positive operators $\Op{H}_{1}$ and
  $\Op{H}_{2}$, resp., (having pure point spectrum) then
  $\Op{H}_{1}$ and $\Op{H}_{2}$ have the same eigenvalues (up to a positive
  constant) and von Neumann
  multiplicities, i.e. there are unitarily equivalent in $\M_{0}$.
\end{lem}

\begin{proof}
  According to \lemref{lem5:0} there is an automorphism $\alpha=\ad\Op{W}$ of
  $\M_{0}$, s.t. $\Op{H}_{2}=\Op{W}\Op{H}_{1}\Op{W}^{*}$. Now the spectrum of
  an operator is invariant under automorphisms, and \lemref{lem2:5} shows that
  also the von Neumann multiplicities are invariant under automorphisms,
  i.e. $\Op{H}_{1}$ and $\Op{H}_{2}$ have the same eigenvalues and von Neumann
  multiplicities, and the unitary equivalence follows from \propref{prop5:2}.
\end{proof} 

The last two lemmas showed that the eigenvalues and multiplicities are actually
the wished complete set of invariants under automorphisms for selfadjoint
operators having pure point spectrum, i.e. the equivalence classes defined by
the equivalence relation of \defnref{defn5:1} can be characterized by
them. 

The only gap left to fill now is the question whether a given decomposition of
the spectrum of the modular operator $\Delta_{0}$ in the sense of
\eqref{eq5:4} and \eqref{eq5:5} gives rise to a corresponding solution of the
inverse problem. This question is answered by the next  
\begin{lem}\label{lem5:4}
  Let $(\mu_{k},m_{k})_{k\in K}$ be a sequence of pairs of positive reals
    $\mu_{k}>0$ and $m_{k}>0$, s.t.
    \begin{subequations}\label{eq5:6}
      \begin{equation}
        m_{k}=l\frac{1}{N}\quad l=1,\dots,N
        \text{ if $\M_{0}$ is type $I_{N}$ ($N\in\N$)},
      \end{equation}
      \begin{equation}
        m_{k}\in\offenl{0}{1}\text{ if $\M_{0}$ is type $II_{1}$},
      \end{equation}
    and
    \begin{equation}
      \sum_{k\in K}m_{k}=1
    \end{equation}
    and 
    \begin{equation}
      \sum_{k\in K}m_{k}\mu_{k}=1
    \end{equation}
  \end{subequations}
  and the relations \eqref{eq5:4} and
  \eqref{eq5:5} are fulfilled. Then there exists a solution
  $\M=\Op{U}\M_{0}\Op{U}^{*}\in NF_{\M_{0}}(\Delta_{0},\Op{J}_{0},u_{0})$,
  s.t.
  $\Op{U}^{*}\Delta_{0}\Op{U}=\Op{H}\Op{J}_{0}\Op{H}^{-1}\Op{J}_{0}$ 
  and $\Op{H}$ has the eigenvalues and multiplicities $(\mu_{k},m_{k})_{k\in
    K}$ (cf. \cite[prop.4.1]{WolII}).
\end{lem}

For the proof we need the following auxiliary results:
\begin{prop}\label{prop5:4}
  If $(m_{k})$ is countable family of positive reals with $\sum m_{k}=1$, then
  there exists in a type $II_{1}$ von Neumann factor $\M$ a family of pairwise
  orthogonal projections $(\Op{E}_{k})$, s.t. $D(\Op{E}_{k})=m_{k}$ for every
  $k$. 
\end{prop}

\begin{proof}
  We construct the $\Op{E}_{k}$ inductively: Since the range of $D_{\M}$ is
  all of $\abge{0}{1}$, if $\M$ is finite, and $\R_{\le0}$, if $\M$ is
  infinite (cf. \cite[8.4.4]{KRII}) there is a projection in $\M$,
  s.t. $D(\Op{E}_{1})=m_{1}$. 

  Suppose now that for $N\in\N$ the $\Op{E}_{k}$ are pairwise orthogonal with
  $D_{\M_{0}}(\Op{E}_{k})=m_{k}$ ($1\leq k<N$). Setting
  $\Op{F}_{N}:=\Op{Id}-\sum_{k=1}^{N}\Op{E}_{k}$ the restricted algebra
  $\Op{F}_{N}\M\Op{F}_{N}$ is again a type $II_{1}$ factor
  (cf. \cite[Ex. 6.9.16]{KRII}) with the dimension function
  \begin{equation*}
    D_{N}(\Op{F}_{n}\Op{E}\Op{F}_{N}):=
    D_{\M_{0}}(\Op{F}_{n}\Op{E}\Op{F}_{N})/D(\Op{F}_{N})
    \quad\forall\Op{F}_{n}\Op{E}\Op{F}_{N}\in\Op{F}_{N}\M\Op{F}_{N},
  \end{equation*}
  where 
  \begin{equation*}
    D_{\M_{0}}(\Op{F}_{N})=D_{\M_{0}}(\Op{Id}-\sum_{k=1}^{N}\Op{E}_{k})
    =1-\sum_{k=1}^{N}D_{\M_{0}}(\Op{E}_{k})\geq m_{N}.
  \end{equation*}
  With the same argument as above there is again a projection
  $\Op{E}_{N}\in\Op{F}_{N}\M\Op{F}_{N}\subset\M$,
  s.t. $D_{N}(\Op{E}_{N})=D(\Op{F}_{N})^{-1}m_{N}\leq1$. Then
  $D_{\M_{0}}(\Op{E}_{N})$ and $\Op{E}_{N}<\Op{F}_{N}\perp\Op{E}_{k}$ ($1\leq
  k<N$).  
\end{proof}

\begin{prop}\label{prop5:2b}
 Let $\Op{H},\tilde\Op{H}$ be two selfadjoint operators with pure
  point spectrum on a Hilbert space $\H_{0}$ which have the same eigenvalues
  and the corresponding eigenspaces have the same dimension. If there is a
  conjugation $\Op{J}_0$ s.t. $\Op{J}_0\Op{H}\Op{J}_0=\Op{H}^{-1}$ and 
  $\Op{J}_0\tilde\Op{H}\Op{J}_0=(\tilde\Op{H})^{-1}$ then there is a
  unitary $\Op{W}\in U(\H_0)$ such that $\Op{W}$ commutes with
  $\Op{J}_0$ and $\tilde\Op{H}=\Op{W}\Op{H}\Op{W}^*$.
\end{prop}

\begin{proof}
  Since $\Op{J}_0\Op{H}\Op{J}_0=\Op{H}^{-1}$ and
  $\Op{J}_0(\tilde\Op{H})\Op{J}_0=(\tilde\Op{H})^{-1}$ we can arrange the
  eigenvalues $\mu_k$ s.t. $\mu_{-k}=\mu_k^{-1}$ and
  $\Op{E}_{-k}=\Op{J}_0\Op{E}_k\Op{J}_0$ and
  $\Op{F}_{-k}=\Op{J}_0\Op{F}_k\Op{J}_0$. Now chose the $\Op{W}_{k}$ for
  $k>=0$ as in the proof of \propref{prop5:2} and set $\Op{\tilde
    W}_k:=\Op{W}_{k}$ for $k>0$,
  $\Op{\tilde W}_k:=\Op{J}_0\Op{W}_{-k}\Op{J}_0$ for $k<0$,  and $\Op{\tilde
    W}_0:=\Op{W}_{0}+\Op{J}_0\Op{W}_{0}\Op{J}_0$ s.t.
  \begin{gather*}
    \Op{\tilde W}_{k}^{*}\Op{\tilde W}_{k}=
    \Op{J}_0\Op{W}_{-k}^{*}\Op{J}_0\Op{J}_0\Op{W}_{-k}\Op{J}_0=
    \Op{J}_0\Op{E}_{-k}\Op{J}_0=\Op{E}_{k}\\
    \text{and}\\
    \Op{\tilde W}_{k}\Op{\tilde W}_{k}^{*}=
    \Op{J}_0\Op{W}_{-k}\Op{J}_0\Op{J}_0\Op{W}_{-k}^{*}\Op{J}_0=
    \Op{J}_0\Op{F}_{-k}\Op{J}_0=\Op{F}_{k}.
  \end{gather*}
  
  Then $\Op{W}:=so-\sum_{k\in K}\Op{W}_k\in U(\H_{0})$ commutes with
  $\Op{J}_0$ and have the stated properties. 
\end{proof}

\begin{proof}[Proof of \lemref{lem5:4}]
  Let $\Op{E}_{k}\in\M_{0}$ ($k\in K$) a family of orthogonal projections in
  $\M_{0}$ with
  $D_{\M_{0}}(\Op{E}_{k})=\tr(\Op{E}_{k})=m_{k}$ and
  $\sum_{k\in K}\Op{E}_{k}=\Op{Id}$  (such a
  family exists according to \propref{prop5:4}). Then we define
  $\Op{H}:=\sum_{k\in K}\mu_{k}\Op{E}_{k}$ which is a non-singular positive
  selfadjoint operator affiliated with $\M_{0}$, has eigenvalues $\mu_{k}$ and
  \begin{equation*}
    \tr(\Op{H})=\sum\mu_k\tr(\Op{E}_k)=\sum\mu_k m_k=1<\infty.
  \end{equation*} 

  Then $(\Delta:=\Op{H}\Op{J}_{0}\Op{H}^{-1}\Op{J}_{0},\Op{J}_0)$ are the 
  modular objects corresponding to $(\M_0,u)$, where 
  $u\in\H_{0}$ is the cyclic and separating vector corresponding to to
  non-singular operator $\Op{T}_{u}:=\Op{H}^{1/2}\Op{V}$ ($\Delta_0$ is the
  modular operator corresponding to $\Op{T}_{u_0}=\Op{H}_0^{1/2}\Op{V}$), and
  $\Delta$ has the same eigenvalues and multiplicities like $\Delta_{0}$ (see
  the proof of \lemref{lem5:1}). This means that they have the same unitary
  invariants in the type $I_\infty$ von Neumann factor $L(\H_0)$

  According to \propref{prop5:2b} there is a unitary $\Op{W}\in U(\H_0)$ s.t.
  $\Delta_{0}:=\Op{W}\Delta\Op{W}^{*}$ and $\Op{W}$ commutes with $\Op{J}_0$ 
  (Since $\Delta$ and $\Delta_0$ both are modular objects with modular 
  conjugation $\Op{J}_0$ we have $\Op{J}_0\Delta_0\Op{J}_0=\Delta_0^{-1}$ and 
  $\Op{J}_0\Delta\Op{J}_0=\Delta^{-1}$). 

  Now we are in exactly the same situation as in the proof of
  \cite[prop.4.1]{WolII} and can show like there that there is also a unitary
  $\Op{U}$ s.t. $\Op{U}$ commutes with $\Op{J}_0$,
  $\Op{U}^*\Delta_0\Op{U}=\Delta$, and $\Op{U}^*u_0=u$, whence
  $\Op{U}\M_0\Op{U}^*$ is a solution with the stated properties.
\end{proof}

Now we can summarize the lemmas of this section in the following
\begin{thm}\label{thm5:1}
  Let $\M_{0}$ be a finite von Neumann factor with cyclic and separating vector
  $u_{0}$ and $\Op{T}_{u_{0}}=\Op{H}_{0}^{-1/2}\Op{V}$ the operator
  corresponding to $u_{0}$. If $\Op{H}_{0}$ has pure point spectrum, also
  $\Delta_{0}$ have it. In this case let $(\lambda_{j})$ $(j\in J)$ be the
  eigenvalues of $\Delta_{0}$. Then  
  \begin{enumerate}
  \item Two solutions $\M_{1},M_{2}\in
    NF_{\M_{0}}(\Delta_{0},\Op{J}_{0},u_{0})$ of the inverse problem with
    corresponding invertible operators $\Op{H}_{i}\eta\M_{0}$ $(i=1,2)$ having
    pure point spectrum are equivalent
    iff $\Op{H}_{1}$ and $\Op{H}_{2}$ have the same eigenvalues and (von
    Neumann) multiplicities.
  \item A positive invertible operator $\Op{H}\eta\M_{0}$ with pure point
    spectrum gives rise to a solution of
    the inverse problem iff its eigenvalues and multiplicities satisfy
    \eqref{eq5:4}, \eqref{eq5:5}, and \eqref{eq5:6}.
  \item When the corresponding operators $\Op{H}$ has pure point spectrum the
    equivalence classes of $\sim$ are completely classified by the spectrum of
    the corresponding operators, i.e. by
    sequences of pairs of positive reals $(\mu_{k},m_{k})$ satisfying
    \eqref{eq5:4}, \eqref{eq5:5}, and \eqref{eq5:6}.
  \end{enumerate}
\end{thm}

\begin{bsp}\label{bsp5:1}
  Here we want to give some examples to illustrate \thmref{thm5:1}.
  \begin{enumerate}
  \item In \cite{WolII} you can find some examples for the type $I$ case.
  \item Let $(\mu_{k},m_{k})$ be the eigenvalues of a positive operator
    $\Op{H}_{0}$ affiliated with a finite factor fulfilling
    the conditions \eqref{eq5:6} s.t. also $(c\mu_{k}^{-1},m_{k})$ in place of
    $(\mu_{k},m_{k})$ fulfill
    conditions \eqref{eq5:6}, where $c>0$ is an appropriate chosen
    constant. Then
    $\Delta_{0}^{-1}=\Op{H}_{0}^{-1}\Op{J}_{0}\Op{H}_{0}\Op{J}_{0}$ is a
    modular operator, the class $NF_{\M_{0}}^{2}(\Delta_{0},\Op{J}_{0},u_{0})$
    exists and is characterized by $(c\mu_{k}^{-1},m_{k})$. Note that, if there
    is a permutation $\sigma$
    s.t. $(c\mu_{\sigma(k)}^{-1},m_{\sigma(k)})=(\mu_{k},m_{k})$ this class is
    just the trivial one, i.e. $NF_{\M_{0}}^{2}=NF_{\M_{0}}^{1}$
    (cf. \corref{cor4:1}).
  \item Let 
    \begin{equation*}
      (10^{-3},10^{-2},10^{-1},1,10,10^{2},10^{3})
    \end{equation*} 
    be the eigenvalues of a modular operator for a type $II_{1}$ factor. Then
    \begin{equation*}
      ((c_{1}\cdot 1,1/4),(c_{1}\cdot 10^{-1},1/4),(c_{1}\cdot
      10^{-2},1/4),(c_{1} \cdot 10^{-3},1/4)),
    \end{equation*}
    \begin{equation*}
      ((c_{2}\cdot 10^{3},1/4),((c_{2}\cdot 10^{2},1/4),(c_{2}\cdot
      10^{1},1/4),((c_{2}\cdot 1,1/4)),
    \end{equation*}
    \begin{equation*}
      ((c_{3}\cdot 1,1/3),(c_{3}\cdot 10^{-1},1/3),(c_{3}\cdot 10^{-3},1/3)),
    \end{equation*}
    and
    \begin{equation*}
      ((c_{4}\cdot 10^{3},1/3),(c_{4}\cdot 10^{1},1/3),(c_{4}\cdot 1,1/3))
    \end{equation*}
    characterize four different classes of solutions of the inverse problem,
    where $c_{i}$ ($i=1,2,3,4$) again are appropriate chosen constants. This
    shows that in this case there are more than the two simple classes of
    solutions of the inverse problem. 
  \item Let $(\mu_{k},m_{k})_{k\in K}$ characterize a class of solutions of
    the inverse problem in the type $II_{1}$ case, where $K$ is a finite index
    set and $m_{l}\not=m_{k}$ for at least one pair $k,l\in K$, then for every
    permutation $\sigma$ of $K$ interchanging $k$ and $l$ also
    $(c\mu_{k},m_{\sigma(k)})$ characterize another class of solutions of
    the inverse problem ($c>0$ a norming constant), where if $\abs{K}=2$ and
    $\mu_{1}=\mu_{2}^{-1}$ this is just the second simple class
    $NF^{2}_{\M_{0}}$, else it is really a new one.
  \item Let again $(\mu_{k},m_{k})_{k\in K}$ be a solution of the inverse
    problem in the type $II_{1}$ case, and let $k,l\in K$ be a pair of
    indices and $\epsilon>0$. Then we get another class by adding $\epsilon$
    to $m_{k}$ and subtracting it from $m_{l}$ where again if $\abs{K}=2$,
    $\mu_{1}=\mu_{2}^{-1}$, and $m_{1}=m_{2}-\epsilon$ this is just the
    second class, else we have really a new one.
  \end{enumerate}
\end{bsp}

\begin{rem}\label{rem5:2}
\begin{enumerate}
  \item \bspref{bsp5:1}.4 and \bspref{bsp5:1}.5 shows that in the type
    $II_{1}$ case if $\Op{H}_{0}$ has more than one eigenvalue, we can always
    construct a third class of solutions, different from the two simple
    classee discussed in \secref{sec:4},
    i.e. $NF_{\M_{0}}\not=NF_{\M_{0}}^{1}\cup NF_{\M_{0}}^{2}$, in contrast to
    the type $I$ case, where for modular operators with generic spectrum we
    have $NF_{\M_{0}}=NF_{\M_{0}}^{1}\cup NF_{\M_{0}}^{2}$ (cf. \cite{WolII}).
  \item
    Unfortunately the classification result presented here applies only to
    operators with pure point spectrum. Whereas in general there are also
    operators with more complicated spectrum (cf. \remref{rem4:0}), for type
    $I$ factors this is no restriction, since all operators generating modular
    operators are trace class operators, hence have pure point spectrum.
  \end{enumerate}
\end{rem}

\begin{appendix}
\section{A Lemma concerning conjugations}
\label{sec3}

For the construction in the next section we need the following lemma
concerning conjugations.

\begin{lem}\label{lem3:1}
Let $\Delta$ be a positive operator on a Hilbert space $\H$, $\Op{J}$ a
conjugation on the same Hilbert space, s.t.
\begin{equation*}
  \Op{J}\Delta\Op{J}=\Delta^{-1}.
\end{equation*}
Further suppose that there are two vectors $v_{1},v_{2}\in\H$, s.t. 
\begin{equation*}
  \Delta v_{i}=\Op{J}v_{i}=v_{i}\quad(i=1,2).
\end{equation*}
Then there is a conjugation $\Op{I}$, s.t.
\begin{equation*}
  \Op{I}\Delta\Op{I}=\Delta,\quad\Op{IJI}=\Op{J}\text{ and }\Op{I}v_{i}=v_{i}
  \quad(i=1,2).
\end{equation*}
\end{lem}

For the proof of this lemma we need some preparatory results.
\begin{prop}\label{prop3:1}
Let $\Delta$ be a selfadjoint operator on a Hilbert space $\H$. Then there is
a conjugation $\Op{K}$, s.t. $\Delta$ is $\Op{K}$-real, i.e.
\begin{equation*}
  \Op{K}\Delta\Op{K}=\Delta
\end{equation*}
(cf. \cite[p.223, ex. 8.1]{Wei}).
\end{prop}

\begin{proof}
  Since $\Delta$ is s.a., there exists a measure space
  $(\Omega,\mathfrak{A},\mu)$, a unitary $\Op{U}:\H\to L_{2}(\mu)$ and a real
  valued 
  measurable function $g$, s.t. for every $f\in L_{2}(\mu)$ with
  $\Op{U}^{*}f\in\Dom{\Delta}$ 
  \begin{equation*}
    (\Op{U}\Delta\Op{U}^{*}f)(t)=g(t)f(t)\text{ $\mu$-a.e. on $\Omega$}.
  \end{equation*}
  Define now
  \begin{equation*}
    (\Op{\tilde K}f)(t)=\kon{f}(t).
  \end{equation*}
  Then $\Op{K}:=\Op{U}^{*}\Op{\tilde K}\Op{U}$ is the wished conjugation.
\end{proof}

\begin{prop}\label{prop3:2}
  Let $\Delta$ be a $\Op{K}$-real s.a. operator on $\H$, where $\Op{K}$ is a
  conjugation. Then:
  \begin{enumerate}
    \item $\mathcal{K}:=\{u\in\H:u=\Op{K}u\}$ is a real subspace of $\H$,
      s.t. $\H=\mathcal{K}+i\mathcal{K}$ and
      $\SProd{\cdot}{\cdot}_{\mathcal{K}}=\SProd{\cdot}{\cdot}$ is a real scalar product on
      $\mathcal{K}$.
    \item $\Delta_{\mathcal{K}}$ is a s.a. operator on the real vector space
      $\mathcal{K}$.
  \end{enumerate}
\end{prop}

\begin{proof}
  \begin{enumerate}
    \item Let $v\in\H$. Then
      \begin{equation*}
        v=\frac{v+\Op{K}v}{2}+i\frac{v-\Op{K}v}{2i}
      \end{equation*}
      and
      \begin{equation*}
        \SProd{u}{v}=\SProd{Ju}{Jv}=\SProd{v}{u}
      \end{equation*}
      for every $u,v\in\mathcal{K}$.
    \item Let $u\in\mathcal{K}\cap\Dom{\Delta}$. Then
      \begin{equation*}
        \Op{J}\Delta u=\Delta\Op{J}u=\Delta u.
      \end{equation*}
      This shows that $\Delta(\mathcal{K}\cap\Dom{\Delta})\subset\mathcal{K}$,
      and the rest follows from standard calculation.
  \end{enumerate}
\end{proof}

\begin{prop}\label{prop3:3}
  Let $\Delta,\Op{J},\H,\raum{K}$ be as in \propref{prop3:2}. Then there is
  an ONB $\{u_{k}\}_{k\in\N}\subset\Dom{\Delta}$ for the real vector space
  $\mathcal{K}$, and $\{u_{k}\}$ is also an ONB for the complex vector space
  $\H$ with $\SProd{\Delta u_{k}}{u_{l}}\in\R$ for $k,l\in\N$.
\end{prop}

\begin{proof}
  This follows immediately from \propref{prop3:2}.
\end{proof}

Now we can prove \lemref{lem3:1}.

\begin{proof}[Proof of \lemref{lem3:1}]
  \begin{enumerate}
    \item Let $\Op{E}(S)$ be the spectral measure of $\Delta$. Then we
      decompose $\H$ in the direct sum of the following three orthogonal
      subspaces:
      \begin{equation*}
        \H=\raum{K}_{-1}\oplus\raum{K}_{0}\oplus\raum{K}_{1},
      \end{equation*}
      where $\raum{K}_{-1}:=\Op{E}(\{\lambda<1\})$,
      $\raum{K}_{0}:=\Op{E}(\{\lambda=1\})$,
      $\raum{K}_{1}:=\Op{E}(\{\lambda>1\})$. Then
      $\Delta(\raum{K}_{j}\cap\Dom{\Delta})\subset\raum{K}_{j}$ ($j=-1,0,1$). From
      $\Op{J}\Delta\Op{J}=\Delta^{-1}$ we see
      $\Op{J}\raum{K}_{j}\subset\raum{K}_{-j}$ ($j=-1,0,1$) and, since
      $\Op{J}^{2}=Id$, also equality holds.
    \item In $\raum{K}_{0}$ we set $\Op{I}_{0}:=\Op{J}_{\raum{K}_{0}}$. Then
      \begin{equation}\label{eq3:1}
        \Op{I}_{0}\Delta_{\raum{K}_{0}}\Op{I}_{0}=\Op{Id}_{\raum{K}_{0}}=
        \Delta_{\raum{K}_{0}}
      \end{equation}
      and
      \begin{equation}\label{eq3:2}
        \Op{I}_{0}\Op{J}_{\raum{K}_{0}}\Op{I}_{0}=\Op{J}_{\raum{K}_{0}}
      \end{equation}
      and, since $v_{i}\in\raum{K}_{0}$,
      \begin{equation}\label{eq3:3}
        \Op{I}_{0}v_{i}=v_{i}\quad(i=1,2).
      \end{equation}
    \item In $\raum{K}_{1}$ we chose according to \propref{prop3:3} an ONB
      $\{u_{k}\}_{k\in\N}\subset\Dom{\Delta}$. Setting
      $u_{-k}:=\Op{J}u_{k}\in\raum{K}_{-1}$ ($k\in\N$) we see that
      $\{u_{-k}\}\subset\Dom{\Delta}$ is an ONB in $\raum{K}_{-1}$. Define now
      the following conjugation $\Op{I}_{j}$ in $\raum{K}_{j}$ ($j=-1,1$)
      \begin{equation}\label{eq3:4}
        \Op{I}_{j}u_{jk}:=u_{jk}.
      \end{equation}
      Then we can calculate with
      $v_{j}=\sum_{l}\alpha_{jl}u_{jl}\in\raum{K}_{j}$:
      \begin{equation*}
        \begin{split}
          \SProd{\Op{I}_{j}\Delta_{\raum{K}_{j}}\Op{I}_{j}u_{jk}}
                {\sum_{l}\alpha_{jl}u_{jl}}&=
          \SProd{\Op{I}_{j}\sum_{l}\alpha_{jl}u_{jl}}
                {\Delta_{\raum{K}_{j}}\Op{I}_{j}u_{jk}}\\
          &=\sum_{l}\kon{\alpha_{jl}}\underbrace{\SProd{u_{jl}}
            {\Delta_{\raum{K}_{j}}u_{jk}}}
                       _{\in\R\text{ (s.\propref{prop3:3})}}\\
          &=\sum_{l}\kon{\alpha_{jl}}\SProd{\Delta_{\raum{K}_{j}}u_{jk}}
                                           {u_{jl}}\\            
          &=\SProd{\Delta_{\raum{K}_{j}}u_{jk}}{v_{j}}.
        \end{split}
      \end{equation*}
      By linear continuation follows for $j=-1,1$
      \begin{equation}\label{eq3:5}
        \Op{I}_{j}\Delta_{\raum{K}_{j}}\Op{I}_{j}=\Delta_{\raum{K}_{j}}.
      \end{equation}
    \item Now we can set
      $\Op{I}:=\Op{I}_{-1}\oplus\Op{I}_{0}\oplus\Op{I}_{1}$. With this
      definition we deduce from \eqref{eq3:1} and \eqref{eq3:5}
      \begin{equation*}
        \Op{I}\Delta\Op{I}=\Delta
      \end{equation*}
      and setting $u=u_{0}+\sum_{j=-1,1;k}\alpha_{jk}u_{jk}\in\H$ from
      \eqref{eq3:2} and \eqref{eq3:4}
      \begin{equation*}
        \begin{split}
          \Op{IJI}u&=\Op{J}u_{0}+
                  \sum_{j,k}\kon{\alpha_{jk}}\Op{I}_{j}\Op{J}\Op{I}_{j}u_{jk}\\
                   &=\Op{J}u_{0}+\sum_{j,k}\kon{\alpha_{jk}}u_{-jk}\\
                   &=\Op{J}u_{0}+\Op{J}\sum_{j,k}\alpha_{jk}u_{jk}=\Op{J}u.
        \end{split}
      \end{equation*}
      Also we see from \eqref{eq3:3}
      \begin{equation*}
        \Op{I}v_{i}=\Op{J}v_{i}=v_{i}\quad(i=1,2)
      \end{equation*}
      and the proof is complete.
  \end{enumerate}
\end{proof}

\section{Some technical results}

\begin{lem}\label{lem2:5}
  Let $\M$ be a finite von Neumann factor with trace $\tr$ and
  $\alpha\in\aut(\M)$ an automorphism. Then
  \begin{equation*}
    \tr(\alpha(\Op{A}))=\tr(\Op{A})\quad\forall\Op{A}\in\M.
  \end{equation*}
\end{lem}

\begin{proof}
  Set $\varphi:=\tr\circ\alpha$. Then an easy calculation shows that also
  $\varphi$ is a trace with
  $\varphi(\Op{Id})=\tr(\alpha(\Op{Id}))=\tr{\Op{Id}}=1$. Since the trace is
  unique we get $\varphi=\tr$, i.e $\tr(\alpha(\Op{A}))=\tr(\Op{A})$ for all
  $\Op{A}\in\M$. 
\end{proof}

\begin{prop}\label{prop4:1}
  Let $\M_{0}$ be a von Neumann factor on the Hilbert space $\H_{0}$. Let
  further $\Delta_{0}=\Op{H}\Op{H}^{'}=\Op{G}\Op{G}^{'}$ a positive
  operator 
  on $\H_{0}$ with $0<\Op{H},\Op{G}\eta\M_{0}$ and
  $0<\Op{H}^{'},\Op{G}^{'}\eta\M_{0}^{'}$. Then $\Op{H}=c\Op{G}$ and
  $\Op{H}^{'}=c^{-1}\Op{G}^{'}$, $c>0$. 
\end{prop}

\begin{proof}
  Since $\Delta_{0}$ is positive, we can examine the unitary group
  $\Delta_{0}^{it}=\Op{H}^{it}\Op{H'}^{it}=\Op{G}^{it}\Op{G'}^{it}$. This
  equality gives
  $\Op{G}^{-it}\Op{H}^{it}=\Op{G'}^{it}\Op{H'}^{-it}\in\alg{C}(\M_{0})=\C$.
  Since also $\Op{G}^{-it}$ and $\Op{H}^{it}$ commute (what is shown by an easy
  computation),  $\Op{G}^{-it}\Op{H}^{it}$ is also a unitary group in $\C$,
  i.e.
  \begin{equation*}
    \Op{G}^{-it}\Op{H}^{it}=c^{it}\text{ with $c>0$}.
  \end{equation*}
  Now the assertion follows from the uniqueness of the generator of a group.
\end{proof}

\end{appendix}

\textbf{Acknowledgements}

The author thanks professor M. Wollenberg for discussing and his usefull hints
and the DFG and the Graduiertenkolleg for financial support.

\bibliographystyle{alpha}

\end{document}